\documentclass[reqno]{amsart}

\usepackage{enumitem}
\usepackage{mathrsfs, color}
\usepackage{amsfonts}
\usepackage{amssymb,amsmath}
\usepackage{amsthm}
\usepackage{graphics}
\usepackage{graphicx}
\usepackage{epsfig}
\usepackage{bm}
\usepackage{subfig}

\usepackage[bookmarks=true]{hyperref}

\allowdisplaybreaks[1]

\newtheorem{theorem}{Theorem}[section]
\newtheorem{lemma}[theorem]{Lemma}

\theoremstyle{definition}
\newtheorem{definition}[theorem]{Definition}

\theoremstyle{remark}
\newtheorem{remark}[theorem]{Remark}

\numberwithin{equation}{section}

\usepackage{algorithm} 
\usepackage{multirow} 
\usepackage{xcolor}
\usepackage{makecell}
\usepackage{lipsum}
\usepackage{amsfonts}
\usepackage{graphicx}
\usepackage{algorithmic}
\ifpdf
  \DeclareGraphicsExtensions{.eps,.pdf,.png,.jpg}
\else
  \DeclareGraphicsExtensions{.eps}
\fi
\renewcommand\arraystretch{1.2}

\begin{document}
\setlength\abovedisplayskip{2pt}
\setlength\belowdisplayskip{1pt}
\setlength\abovedisplayshortskip{2pt}
\setlength\belowdisplayshortskip{1pt}

\title[Infinite-Dimensional Continuous Normalizing Flow]
 {Learning Informative Prior with Infinite-Dimensional Continuous Normalizing Flow for Bayesian Inverse Problem}

\author[Y. Zhao]{Yang Zhao}
\address{School of Mathematics and Statistics, Xi'an Jiaotong University,
	Xi'an, Shaanxi 710049, China}
\email{zhaoyangxt@xjtu.edu.cn}

\author[J. Jia]{Junxiong Jia}
\address{School of Mathematics and Statistics,
Xi'an Jiaotong University,
 Xi'an
710049, China}
\email{jjx323@xjtu.edu.cn}

\author[T. Zhou]{Tao Zhou}
\address{Academy of Mathematics and Systems Sciences, Chinese Academy of Sciences,
    Beijing, 100190, China}
\email{tzhou@lsec.cc.ac.cn}

\subjclass[2020]{35R30, 65N21, 62F15}

\date{}

\keywords{functional continuous normalizing flow, inverse problems for partial differential equation, infinite-dimensional Bayesian analysis, prior learning}

\begin{abstract}
Infinite-dimensional Bayesian inference methods, as a core class of methodologies for solving inverse problems of partial differential equations, have been widely applied across numerous scientific fields. A critical obstacle to constructing effective Bayesian inversion methods is how to construct an informative prior measure, particularly in the infinite-dimensional setting. Based on an infinite-dimensionally well-defined neural ordinary differential equation, we establish a general learning framework for infinite-dimensional nonlinear transformations, which can be incorporated into two practical Bayesian prior-learning paradigms corresponding to the direct and indirect data scenarios. A key ingredient of this framework is an equivalence between non-Gaussian and Gaussian measures, which is highly nontrivial: nonlinear transformations with genuinely infinite-dimensional range are involved, and several new limiting mechanisms have to be clarified. The developed framework is applied to three representative inverse problems -- simple source inversion, inverse obstacle scattering, and inverse heat conduction. All three examples substantiate the effectiveness of the proposed methods and, at the same time, demonstrate their broad applicability.
\end{abstract}

\maketitle



\section{Introduction}\label{sec:introduction}
Estimating infinite-dimensional parameters for inverse problems governed by partial differential equations (PDEs) from measurement data is a fundamental problem in various scientific and engineering fields \cite{Kirsch2011Book,Benning2018ActaNum}. 
In recent years, there has been growing emphasis on statistical approaches that not only provide approximate solutions but also quantify the associated uncertainties, a capability crucial for tasks such as artifact detection \cite{Zhou2020SJIS}. 

Unlike deterministic approaches, Bayesian inference offers a unified framework to quantify uncertainty in inverse problems \cite{bishop2006pattern,kaipio2006statistical,dashti2013bayesian}, originally developed in finite-dimensional spaces \cite{bishop2006pattern,kaipio2006statistical}. However, finite-dimensional methods often fail to characterize the continuous nature of inverse problems of PDEs. To retain the inherent infinite-dimensional structure of PDE-based inverse problems, defining Bayesian inference directly in infinite-dimensional spaces has been widely advocated \cite{jia2021variational,cotter2013,Jia2022VINet,sui2024noncentered,pinski2015algorithms,Pinski2015SIAMMA,dashti2013bayesian}. This paper concerns infinite-dimensional Bayesian inversion \cite{dashti2013bayesian,stuart2010inverse}. Two central challenges in this field are the construction of informative prior measures and the development of efficient numerical algorithms, both essential for accurate and feasible inversion in infinite dimensions. In this work, we focus on the construction of prior measures.

Effective priors are essential for algorithmic efficiency, as they act as regularizers and encode prior information about model parameters. A common strategy for constructing structured priors  uses convergent series expansions with coefficients sampled from finite-dimensional distributions, yielding priors such as uniform and Gaussian priors \cite{dashti2013bayesian,stuart2010inverse}. Another popular choice is the edge-preserving Besov prior \cite{MR2558305}, often employs Haar wavelet bases, which can impose restrictive assumptions \cite{arridge2019solving}. As an alternative, Markkanen et al. proposed the Cauchy difference prior \cite{cauchydiff}, which can better adapts to complex geometric structures in inverse problems. Although these priors characterize the smoothness of function-valued parameters, they often fail to fully capture the problem-specific structures and the historical information inherent in certain inverse problems.

To overcome the limitations of handcrafted priors, data-driven prior-learning schemes have been proposed, most of which demand explicit parameter samples in Euclidean settings. For example, Asim et al. employ normalizing flow based invertible networks as learned image priors to capture realistic signals and mitigate reconstruction bias for denoising, compressive sensing and inpainting \cite{asim2020invertible}. Cai et al. propose NF-ULA, a normalizing flow-driven unadjusted Langevin algorithm that learns explicit image priors for stable Bayesian sampling with theoretical convergence guarantees \cite{MR4728740}. These generative learning methods excel at modeling complex distributions and have become widely adopted replacements for traditional handcrafted priors. Nevertheless, these methods rely on an experimental setup where full training samples of underlying model parameters can be collected. This favorable data condition is often difficult to achieve in practical PDE inverse problems, where model parameters cannot be measured directly and only indirect observations through forward operators are available.

To tackle this problem, many studies have focused on indirect prior learning, which extracts shared prior knowledge from measurement data across related inverse problems governed by the same PDE \cite{daras2023ambient,gao2023image,akyildiz2025efficient}. For example, Akyildiz et al. \cite{akyildiz2025efficient} developed a scalable framework to learn priors from noisy indirect measurements. It leverages generative models and the sliced-Wasserstein distance for distribution matching, and employs neural operators to reduce the heavy computational cost of PDE simulations.  Jia et al. \cite{jia2026nonparametric} proposed a nonparametric Bayesian prior learning approach based on meta-learning and PAC-Bayesian theory. It constructs data-adaptive priors directly from observations and establishes generalization bounds for both linear and nonlinear PDE systems, leading to improved efficiency in posterior sampling and more reliable uncertainty quantification.

Beyond the design of training schemes with direct and indirect dataset, constructing well-posed prior distributions over infinite-dimensional spaces is another critical task. Since the Lebesgue  measure is not well defined on infinite-dimensional spaces, conventional probability density functions--interpreted as the Radon Nikodym (RN) derivative of a target measure with respect to Lebesgue measure--no longer exist in such settings. However, density functions are indispensable for both posterior sampling (e.g., Markov chain Monte Carlo algorithms) and prior learning (e.g., learning the prior via negative log-likelihood).  This motivates us to seek a valid surrogate for Lebesgue-type densities on infinite-dimensional spaces.  A standard approach is to introduce a common reference measure $\mu_0$, with respect to which our target $\mu$ is absolutely continuous, allowing $\mu$ to be formulated via the RN derivative \cite{dashti2013bayesian}. Notably, infinite-dimensional measures frequently exhibit mutual singularity \cite{gaussianmeasure}, breaking absolute continuity. Hence, constructing valid infinite-dimensional priors becomes challenging, as they must be absolutely continuous with respect to the reference measure.


In this paper, we propose a novel modeling framework: infinite-dimensional functional continuous normalizing flows (f-CNFs) for Bayesian inversion. This framework is tailored to the construction of Bayesian priors. By learning from historical datasets of the PDE model, the generated prior (hereafter referred to as the f-CNF prior) enables the Bayesian framework to yield a more accurate posterior, particularly when the likelihood is uninformative. Our f-CNF architecture is rigorously constructed to transform a simple reference measure into a complex target measure (serving as the prior) via a parameterized neural ordinary differential equation (ODE). This construction preserves measure equivalence, allowing the resulting prior to be characterized via the RN derivative. By optimizing the f-CNF parameters, we can incorporate historical information into the Bayesian prior.

In our numerical experiments, we design the learning of infinite-dimensional priors for two data types: direct parameter samples and indirect observational data. Learning from direct data constructs an infinite-dimensional generative model, which subsequently serves as a Bayesian prior. Consequently, unlike existing infinite-dimensional generative models \cite{diffusionf2, kerrigan2023functional, shikhman2026discretization}, we must guarantee absolute continuity and ensure that the corresponding RN derivative is analytically tractable; see Subsection \ref{subsec:RN derivative} for details. For indirect-data training, we adopt a Wasserstein loss similar to \cite{akyildiz2025efficient}. While their work optimizes a limited set of parameters, we train the full f-CNF network to achieve improved flexibility. Upon completing the learning of the f-CNF prior, we combine it with the preconditioned Crank--Nicolson (pCN)  \cite{cotter2013} and the Sequential Monte Carlo with Gaussian Mixture (SMC-GM)  \cite{lu2025sequential} algorithms to construct dedicated posterior samplers. Three numerical experiments verify the efficacy of our prior-learning and sampling schemes.

In summary, this work mainly contains three contributions:
\begin{enumerate}
    \item To address the challenge of learning priors over functions outlined in \cite{arridge2019solving}, we establish a general measure equivalence theorem for the nonlinear transformation induced by an infinite-dimensional neural ODE. This theorem leads to a new model named f-CNF, which has the ability to encode intricate information contained in historical datasets.
    \item In contrast to conventional log-likelihood-based approaches \cite{MR4728740,NODE}, we establish a training loss for the proposed f-CNF model that is well-defined in infinite dimensions, ensuring dimension invariance---a crucial property for PDE inverse problems---from the ground up.
    \item For the non-Gaussian prior learned by the f-CNF model, we prove the well-posedness of the posterior and develop efficient posterior sampling algorithms by adapting pCN \cite{cotter2013} and SMC-GM \cite{lu2025sequential}. Numerical experiments on three representative PDE inverse problems demonstrate the effectiveness and generality of the proposed framework.
\end{enumerate}

The remainder of this paper is organized as follows. Section \ref{sec:f-CNF} introduces the f-CNF framework. Section \ref{sec:IP_and_Prior} reviews the infinite-dimensional Bayesian inverse problem, specifies the f-CNF prior, constructs tailored samplers, and develops f-CNF training schemes. Section \ref{sec:wellposedness} analyzes posterior well-posedness with the f-CNF prior. Section \ref{sec:numerical experiment} validates our algorithms via three canonical PDE inverse problems. Section \ref{sec:Conclusion} summarizes the work and outlines future research directions.

\section{Functional Continuous Normalizing Flows}\label{sec:f-CNF}
In this section, we establish the theoretical framework of f-CNF in infinite-dimensional space. To
circumvent the obstacle posed by the singularity of measures in infinite-dimensional space, and enable the ‘density estimation’, we present a novel theoretical framework for the proposed f-CNF model and provide explicit computations of the corresponding RN derivative.
\subsection{Basic Setting}\label{subsec:Basic setting}
The generative model for data $\mathbf{x} \in \mathbb{R}^D$ proposed by Chen et al., known as continuous normalizing flows \cite{NODE}, is formulated through continuous-time dynamics. 
This concept readily extends to infinite-dimensional spaces. We consider the following ODE in infinite-dimensional space \cite{zeidler2013nonlinear}:
\begin{align}\label{equ:forward}
\begin{cases}
\displaystyle \frac{dv(t)}{dt} = h(v(t),t; \theta), \\
v(0) = u_0
\end{cases}
\end{align}
where $u_0$ denotes an initial function sampled from a Gaussian measure $\mu_0 = \mathcal{N}(0, \mathcal{C}_0)$ on a separable Hilbert space $\mathcal{H}_u$, and $h$ denotes a mapping from $\mathcal{H}_u \times \mathbb{R}$ to $\mathcal{H}_u$ equipped with parameter $\theta \in \Theta$ (here $\Theta$ denotes the space of all trainable neural network parameters). The ODE \eqref{equ:forward} then propagates the initial function \( u_0 \in \mathcal{H}_u \) forward in time, producing the terminal solution \(v(T) = u_T \in \mathcal{H}_u \) at time \( T \). For a fixed parameter $\theta\in\Theta$, suppose that $h(u,t;\theta)$ is continuous on $\mathcal{H}_u\times[0,T]$, globally Lipschitz continuous with respect to $u$ uniformly in $t$ (i.e., $\|h(u_1,t;\theta)-h(u_2,t;\theta)\|_{\mathcal{H}_u}\leq L\|u_1-u_2\|_{\mathcal{H}_u}$), and uniformly bounded (i.e., $\|h(u,t;\theta)\|_{\mathcal{H}_u}\leq K$), where constants $L$ and $K$ are independent of $u$ and $t$. Then, the ODE \eqref{equ:forward} admits a unique solution on $[0,T]$ \cite{zeidler2013nonlinear}. In this situation, the entire ODE propagates process, transforming $u_0$ to $u_T$, can be regard as a mapping, denoted $f_{\theta}$, from the infinite-dimensional Hilbert space $\mathcal{H}_u$ to itself:
\begin{equation}\label{equ:ftheta}
f_{\theta}(u_0) = u_T,
\end{equation}
and $u_T$ can be regard as a sample from the measure $\mu_{f_{\theta}}={f_{\theta}}_{\#}\mu_0$. Here ${f_{\theta}}_{\#}\mu_0$ denotes the pushforward measure $f_\theta$ induced by $\mu_0$ \cite{gaussianmeasure}. For any Borel set $\mathcal{B}$, the pushforward measure is rigorously defined as $({f_{\theta}}_{\#}\mu_0)(\mathcal{B}) := \mu_0\big(f_\theta^{-1}(\mathcal{B})\big)$, where $f_\theta^{-1}(\mathcal{B}) = \{u \mid f_\theta(u) \in \mathcal{B}\}$ stands for the preimage of $\mathcal{B}$ under $f_\theta$.
 Furthermore, $f_\theta$ admits an inverse mapping governed by a backward ODE \cite{zeidler2013nonlinear}, defined via
\begin{align}\label{equ:inverse}
\left\{
\begin{array}{l}
\displaystyle \frac{dv(t)}{dt} = h(v(t),t; \theta), \\
v(T) = u_T.
\end{array}
\right.
\end{align}
This backward flow maps $u_T$ to $u_0=v(0)$, denoted $g_{\theta}(u_T)=u_0$. It is straightforward to verify that $g_{\theta}$ and $f_{\theta}$ are mutual inverses, since the solutions of the ODEs \eqref{equ:forward} and \eqref{equ:inverse} exist and are unique. 

Via the mapping $f_{\theta}$ on the infinite-dimensional Hilbert space $\mathcal{H}_u$, we construct a parameterized measure $\mu_{f_\theta}$ over $\mathcal{H}_u$. Varying $\theta\in\Theta$ endows $\mu_{f_\theta}$ with rich approximation capacity. Our core strategy is to embed historical information of PDE inverse problems into $\mu_{f_\theta}$ and employ $\mu_{f_\theta}$ as the Bayesian prior to enhance posterior accuracy. Before elaborating on this concept, we must first address a fundamental issue concerning infinite-dimensional spaces in the following subsection.

\subsection{The Radon-Nikodym Derivative}\label{subsec:RN derivative}
Due to its infinite-dimensional nature, an infinite-dimensional Hilbert space does not admit a translation-invariant Lebesgue measure. Consequently, the Lebesgue density functions commonly utilized in Euclidean space do not exist in infinite-dimensional space $\mathcal{H}_u$. As a substitute, RN derivative with respect to another common measure is frequently employed in place of the Lebesgue density function to characterize a measure within an infinite-dimensional space \cite{dashti2013bayesian}. Specifically, suppose we wish to characterize a target measure $\mu_{\text{target}}$. A common approach involves first identifying a reference measure $\mu_{\text{simple}}$, ensuring that $\mu_{\text{target}}$ is absolutely continuous with respect to $\mu_{\text{simple}}$ (or satisfying the stronger condition that $\mu_{\text{target}}$ and $\mu_{\text{simple}}$ are equivalent). Under these conditions, $\mu_{\text{target}}$ can be described via the RN derivative $\frac{d\mu_{\text{target}}}{d\mu_{\text{simple}}}(u)$.

It is noteworthy that the RN derivative is meaningful only if the measure $\mu_{\text{target}}$ is absolutely continuous with respect to the measure $\mu_{\text{simple}}$. Considering the push-forward measure $\mu_{f_{\theta}} = {f_{\theta}}_{\#}\mu_0$ constructed in Subsection \ref{subsec:Basic setting}, characterizing it via the RN derivative requires identifying a reference measure with the absolutely continuous property. A natural candidate is the pre-transformed measure $\mu_0$ of our f-CNF model $f_{\theta}$. Consequently, a pertinent question arises: can we design the ODE flow $f_{\theta}$ such that $\mu_{f_{\theta}}$ is absolutely continuous with respect to $\mu_0$?

In contrast to measures in Euclidean spaces, measures in infinite-dimensional spaces are inherently prone to mutual singularity. For $m\notin\mathcal{H}(\mu_0)$ (here $\mathcal{H}(\mu_0)$ denotes the Cameron-Martin space of $\mu_0$ \cite{gaussianmeasure}), the shifted pushforward $\mu_m={f_{m}}_{\#}\mu_0$ with $f_m(u)=u+m$ is singular to $\mu_0$. Likewise, the scaled pushforward $\mu_2={f_{2}}_{\#}\mu_0$ for $f_2(u)=2u$ also yields a singular measure relative to $\mu_0$. These examples illustrate that shifts and scalings generally produce singular Gaussian measures on infinite-dimensional spaces \cite{gaussianmeasure}.
Although simple, the transformations $f_{m}$ and $f_{2}$ map $\mu_0$ to singular measures, violating the absolute continuity required for the RN derivative. Consequently, constructing expressive flows that preserve absolute continuity is nontrivial. We resolve this issue with the general theorem below, which guarantees that $\mu_{f_{\theta}}$ is equivalent to $\mu_0$.

Since continuous normalizing flows can be regarded as the continuous-time limit of discrete normalizing flows \cite{NODE}, we first establish measure equivalence for functional normalizing flows, and then extend it to the continuous f-CNF via a limiting argument. Proofs are given in the Appendix \ref{subappendix:proof of fnf}.


First, we establish conditions ensuring the equivalence of measures before and after transformation for functional normalizing flows. This result refines our previous conclusions in \cite{fnf}, as detailed below.

\begin{lemma}\label{thm:discrete}
Let $\mathcal{H}_u$ be a separable Hilbert space equipped with the Gaussian measure $\mu_0 = \mathcal{N}(0,\mathcal{C}_0)$, and $\mathcal{H}=\mathcal{H}(\mu_0)$ be the Cameron-Martin space of the measure $\mu_0$. For each $n = 1, 2, \ldots, N$, let $f_{\theta_n}^{(n)}(u)=u+\mathcal{F}_{\theta_n}^{(n)}(u)$, where $\mathcal{F}_{\theta_n}^{(n)}$ is an operator from $\mathcal{H}_u$ to $\mathcal{H}_u$. The composite transformation $f_{\theta}$ is given by $f_{\theta}(u)=f_{\theta_N}^{(N)} \circ f_{\theta_{N-1}}^{(N-1)} \circ \cdots \circ f_{\theta_1}^{(1)}(u)$, and $\mu_{f_{\theta}}$ denotes its corresponding push-forward measure by $\mu_{f_{\theta}}={f_{\theta}}_{\#}\mu_0$. For any $\theta \in \Theta$, assume that for all $n = 1,2,\ldots,N$ the following conditions hold:
\begin{itemize}
\item The space $\operatorname{Im}(\mathcal{F}_{\theta_n}^{(n)}) \subset \mathcal{H}$, and $\mathcal{F}_{\theta_n}^{(n)}:\mathcal{H}_u\to\mathcal H$ is continuous.
\item For any $u \in \mathcal{H}_u$, $D_{\mathcal H}\mathcal{F}_{\theta_n}^{(n)}(u)$ is a trace class operator on $\mathcal H$ and continuous in $u$ with respect to the trace norm.
\item The operator $f_{\theta_n}^{(n)}$ is a homeomorphism.
\item For any $u \in \mathcal{H}_u$, all point spectrum of $D_{\mathcal H}\mathcal{F}_{\theta_n}^{(n)}(u)$ are not in $(-\infty,-1]$.
\end{itemize}
Then $\mu_{f_\theta}$ is equivalent to $\mu_0$, with the RN derivative given by
\begin{align*}
    \frac{d\mu_{f_{\theta}}}{d\mu_0}(f_{\theta}(u))
    =\prod\limits_{n=1}^{N}\left | {\det}_{1}(D_{\mathcal H}f_{\theta_n}^{(n)}(u_{n-1})) \right |^{-1}
    \exp\left(\frac{1}{2}\langle \mathcal{T}_{\theta}(u),\mathcal{T}_{\theta}(u)\rangle_{\mathcal{H}}
    +\langle u,\mathcal{T}_{\theta}(u)\rangle_{\mathcal{H}}^{\mathrm{ext}}\right),
\end{align*}
for $\mu_0$-almost every $u\in\mathcal H_u$, where $u_n=f_{\theta_n}^{(n)} \circ f_{\theta_{n-1}}^{(n-1)} \circ \cdots \circ f_{\theta_1}^{(1)}(u)$ for each $n = 1, 2, \ldots, N$, $\mathcal{T}_{\theta}(u)=f_{\theta}(u)-u$， ${\det}_{1}(D_{\mathcal H}f_{\theta_n}^{(n)}(u_{n-1}))$ is the Fredholm determinant of the linear operator $D_{\mathcal H}f_{\theta_n}^{(n)}(u_{n-1})$, and $\langle u,\mathcal T_\theta(u)\rangle_{\mathcal H}^{\mathrm{ext}}:=\operatorname{Tr}_{\mathcal H}
\bigl(D_{\mathcal H}\mathcal T_\theta(u)\bigr)
-\delta_{\mu_0}\mathcal T_\theta(u)$ with $\delta_{\mu_0}$ be the Gaussian divergence.
\end{lemma}
{\color{black}
\begin{remark}
The Fredholm determinant $\det_1(\cdot)$ extends the classical determinant from finite-dimensional spaces to trace-class operators on infinite-dimensional spaces. For the reader's convenience, a rigorous definition is provided in Appendix~\ref{subappendix:Fredholm and Regularized Determinants}.
\end{remark}
}

\begin{remark}\label{remark:ext}
The extended Cameron--Martin pairing
$\langle u,\mathcal T_\theta(u)\rangle_{\mathcal H}^{\mathrm{ext}}$
is defined as the limit of
$\langle \Pi_m u,\Pi_m\mathcal T_\theta(u)\rangle_{\mathcal H}$
as $m\to\infty$, with convergence in $\mu_0$-measure.
Here $\Pi_m$ denotes the orthogonal projection in $\mathcal H_u$
onto the subspace spanned by the first $m$ elements of an
orthonormal eigenbasis of $\mathcal C_0$.
This pairing admits the representation
$
\langle u,\mathcal T_\theta(u)\rangle_{\mathcal H}^{\mathrm{ext}}
=
\operatorname{Tr}_{\mathcal H}
\bigl(D_{\mathcal H}\mathcal T_\theta(u)\bigr)
-\delta_{\mu_0}\mathcal T_\theta(u),
$
where $\delta_{\mu_0}$ denotes the Gaussian divergence
\cite{gaussianmeasure}, and $D_{\mathcal H}$ denotes
differentiation along $\mathcal H$.
\end{remark}

\begin{remark}
When $u \in \mathcal H$, the extended pairing agrees with the ordinary Cameron--Martin inner product. Writing $\mathcal C_0e_k=\lambda_ke_k$, we set $h_k=\sqrt{\lambda_k}\,e_k$, so that $\{h_k\}_{k=1}^{\infty}$ forms an orthonormal basis of $\mathcal H$. For $u\in\mathcal H$, by Theorem 5.8.3 of \cite{gaussianmeasure}, we can choose the following representative of the divergence:
\begin{align*}
\delta_{\mu_0}\mathcal T_\theta(u)
&=
\sum_{k=1}^{\infty}
\langle D_{\mathcal H}\mathcal T_\theta(u)h_k,h_k\rangle_{\mathcal H}-
\sum_{k=1}^{\infty}
\langle u,h_k\rangle_{\mathcal H}
\langle\mathcal T_\theta(u),h_k\rangle_{\mathcal H}\\
&=
\operatorname{Tr}_{\mathcal H}
\bigl(D_{\mathcal H}\mathcal T_\theta(u)\bigr)
-\langle u,\mathcal T_\theta(u)\rangle_{\mathcal H}.
\end{align*}
With this choice, we have $\langle u,\mathcal T_\theta(u)\rangle_{\mathcal H}^{\mathrm{ext}}
=\operatorname{Tr}_{\mathcal H}
\bigl(D_{\mathcal H}\mathcal T_\theta(u)\bigr)
-\delta_{\mu_0}\mathcal T_\theta(u)=\langle u,\mathcal T_\theta(u)\rangle_{\mathcal H}.$
\end{remark}

The f-CNF represents the continuous-time limit of its discrete counterparts. Discretizing the ODE in \eqref{equ:forward} using an Euler scheme recovers a standard functional normalizing flow; allowing the step size to vanish refines this discrete flow into the f-CNF. Below, we develop the accompanying measure-theoretic limit theory, which guarantees that the desirable properties of normalizing flows in infinite-dimensional space persist throughout this limiting procedure.
\begin{theorem}\label{thm:limits}
Let $\mu_0$ be a Gaussian measure on a separable Hilbert space $\mathcal{H}_u$, and let $\{f_N\}_{N \in \mathbb{N}^+}$ be a sequence of measurable mappings from $\mathcal{H}_u$ to itself.
Suppose the following conditions hold:
\begin{itemize}
    \item There exists a mapping $f: \mathcal{H}_u \to \mathcal{H}_u$ such that
          $\lim\limits_{N\to\infty} \|f_N(u) - f(u)\|_{\mathcal{H}_u} = 0$ for any $u \in \mathcal{H}_u$.
    \item Measure $\mu_0$ is equivalent to ${f_{N}}_{\#}\mu_0$ for $N \in \mathbb{N}^+$, with RN derivative $\frac{d{f_{N}}_{\#}\mu_0}{d\mu_0}(u)=\mathcal{R}_N(u)$.
    \item There exists $g\in L^1(\mu_0)$ such that
$|\mathcal{R}_N(u)|\le g(u)$ $\mu_0$-a.s. for all $N$.
    \item The limit $\mathcal{R}(u)=\lim_{N\to\infty}\mathcal{R}_N(u)$
exists and is strictly positive $\mu_0$-a.s.

\end{itemize}
Then $\mu_0$ is equivalent to ${f}_{\#}\mu_0$, with RN derivative given by
$
\frac{d{f}_{\#}\mu_0}{d\mu_0}(u) = \mathcal{R}(u).
$
\end{theorem}

Finally, combining Lemma \ref{thm:discrete} and Theorem \ref{thm:limits}, we derive conditions for continuous normalizing flows to preserve measure equivalence in infinite-dimensional spaces. The necessary function-space definitions precede the main theorem.

\begin{definition}(See \cite{gaussianmeasure} for detail)
Let $\mathcal{H}$ be the Cameron-Martin space of $\mu_0$. Denote by $\mathcal{HC}^1(\mu_0, \mathcal{H})$ the class of all $\mu_0$-measurable mappings $\mathcal{F}: \mathcal{H}_u \to \mathcal{H}$ such that for almost every $u \in \mathcal{H}_u$, the mapping $h \mapsto \mathcal{F}(u + h)$ is Fréchet differentiable along $\mathcal{H}$, and its derivative is a Hilbert--Schmidt operator such that the mapping $h \mapsto D_{\mathcal H}\mathcal{F}(u + h)$ from $\mathcal{H}$ to $\mathfrak{S}_2(\mathcal{H})$ is continuous. Here, $D_{\mathcal H}\mathcal{F}(u)$ denotes the derivative, and $\mathfrak{S}_2(\mathcal{H})$ denotes the Hilbert space of Hilbert--Schmidt operators on $\mathcal{H}$, equipped with the Hilbert--Schmidt norm.
\end{definition}

\begin{definition}
Let $\mu_0$ be the Gaussian measure $\mathcal{N}(0,\mathcal{C}_0)$ on $\mathcal{H}_u$, with injective covariance operator $\mathcal{C}_0$, and let $\mathcal{H}$ denote its Cameron--Martin space. We define $\mathcal{T}(\mu_0,[0,T], \mathcal{H})$ as the collection of functions $h(u,t)$ mapping $\mathcal{H}_u \times [0,T]$ into $\mathcal{H}_u$ such that the following conditions hold:

\begin{itemize}
\vskip 0.1 cm

\item $h$ is continuous on $\mathcal{H}_u\times[0,T]$, and $h(\cdot,t)\in\mathcal{HC}^{1}(\mu_0,\mathcal{H})$ for each $t\in[0,T]$;

\vskip 0.1 cm
\item $\sup\limits_{t \in [0,T],\, u \in \mathcal{H}_u} \left\| \mathcal{C}_0^{-1} h(u,t) \right\|_{\mathcal{H}_u} \leq M_1 < \infty$;
\vskip 0.1 cm

\item $D_{\mathcal H}h(u,t)$ is a trace class operator on $\mathcal H$ and continuous in $(u,t)$ with respect to the trace norm, satisfying
$
\sup_{u,t}\|D_{\mathcal H}h(u,t)\|_1\leq M_2<\infty.
$

\vskip 0.1 cm
\item For each $t \in [0,T]$, $h(\cdot,t)$ is Lipschitz continuous, i.e.,
$
\left\| h(u_1,t) - h(u_2,t) \right\|_{\mathcal{H}_u}
\leq L \left\| u_1 - u_2 \right\|_{\mathcal{H}_u}.
$
\end{itemize}
Here, the constants $M_1$, $M_2$, and $L$ are independent of $u$ and $t$. The notation $\|\cdot\|_1=\operatorname{Tr}(|\cdot|)$ denotes the trace norm on a Hilbert space $\mathcal{H}$ \cite{reed2012methods}.
\end{definition}

Subsequently, we present the measure equivalence theorem for the f-CNF model and derive the corresponding expression for the RN derivative.

\begin{theorem}\label{thm:RN derivative}
Consider the following ODE:
\begin{align}\label{equ:forward2}
\left\{
\begin{array}{l}
\displaystyle \frac{dv(t)}{dt} = h(v(t),t; \theta), \\
v(0) = u.
\end{array}
\right.
\end{align}
Suppose that, for fixed $\theta\in\Theta$,
$h(\cdot,\cdot;\theta)\in
\mathcal{T}(\mu_0,[0,T],\mathcal{H})$.
Let $f_\theta$ denote the flow map defined
analogously to \eqref{equ:ftheta}.
Then, ${f_\theta}_{\#}\mu_0$ is equivalent
to $\mu_0$. Letting \(\mathcal{T}_\theta(u) = f_\theta(u) - u\), the RN derivative is given by
\begin{align*}
\frac{d{f_{\theta}}_{\#}\mu_0}{d\mu_0}(f_{\theta}(u))
&=\exp\left(\frac12\lVert\mathcal{T}_{\theta}(u)\rVert_{\mathcal{H}}^2 + \langle u,\mathcal{T}_{\theta}(u) \rangle_{\mathcal{H}}^{\mathrm{ext}}
-\int_0^{T} \mathrm{Tr}\, D_{\mathcal{H}}h(v(t),t; \theta) dt\right),
\end{align*}
where \(v(t)\) denotes the solution to \eqref{equ:forward2} with initial condition \(v(0)=u\).
Here, the extended inner product
$\langle u,\mathcal{T}_{\theta}(u)\rangle_{\mathcal{H}}^{\mathrm{ext}}$
is defined as in Lemma \ref{thm:discrete}. 
\end{theorem}

Under the conditions of Theorem~\ref{thm:RN derivative}, the ODE \eqref{equ:forward2} is well-posed, and the transformation $f_{\theta}$ pushes the prior measure $\mu_0$ forward to a measure $\mu_{f_{\theta}} = {f_{\theta}}_{\#}\mu_0$ that is equivalent to $\mu_0$. Moreover, its RN derivative can be computed explicitly. In the remainder of this paper, we demonstrate how the proposed f-CNF can be applied to a variety of PDE inverse problems.

\section{Inverse Problem and Bayesian Prior}\label{sec:IP_and_Prior}
\subsection{Bayesian Approach for Inverse Problem}\label{subsec:Bayesian Inverse Problem}

For a large class of inverse problems involving PDEs, sparse measurement data are practically adopted \cite{Cotter2009IP, nickl2020bernstein}, as such data are more easily acquired.  In this subsection, we introduce the foundational concepts of infinite-dimensional Bayesian inverse problems with finite-dimensional sparse data.

Let $\mathcal{H}_u$ and $\mathcal{H}_w$ be separable Hilbert spaces representing the parameter space and the solution space, respectively, and let $N_d$ be a positive integer. Since the measurement is often perturbed by noise, we consider the observation model
\begin{align}\label{eq}
\bm{d} = \mathcal{S}\mathcal{G}(u) + \bm{\epsilon},
\end{align}
where $\bm{d}\in \mathbb{R}^{N_d}$ denotes the measurement data, $u\in \mathcal{H}_u$ is the parameter of interest, $\mathcal{G}$ is the PDE solution operator from $\mathcal{H}_u$ to $\mathcal{H}_{w}$, $\mathcal{S}$ is the measurement operator from $\mathcal{H}_{w}$ to $\mathbb{R}^{N_d}$, and $\bm{\epsilon}$ is Gaussian noise. Based on the framework of infinite-dimensional Bayesian inference \cite{dashti2013bayesian}, we are able to preserve the fundamental Bayes' formula for the inverse problem:
\begin{align}\label{equ:bayes infinite}
	\qquad\,\,\frac{d\mu}{d\mu_{\textbf{prior}}}(u) = \frac{1}{Z_{\mu}}\mathbb{L}(\bm{d}|\cdot),
\end{align}
where $\mathbb{L}(\bm{d}|\cdot):=\exp (-\Phi(\bm{d}|\cdot) ) : \mathcal{H}_u \rightarrow \mathbb{R}$ is the likelihood function, and $Z_{\mu}$ is a positive finite constant given by $Z_{\mu} = \int_{\mathcal{H}_u} \exp (-\Phi(\bm{d}|u))\mu_{\textbf{prior}}(du)$.

The selection of an appropriate prior distribution $\mu_{\textbf{prior}}$ is essential for obtaining an appropriate Bayesian posterior. In this section, we propose a prior based on the f-CNF model introduced in Subsection \ref{subsec:Basic setting}. This prior can incorporate information from historical data, enabling the Bayesian framework to yield a more accurate posterior, particularly when the likelihood function provides limited information.

\subsection{Functional Continuous Normalizing Flows as Prior}\label{subsec:f-CNF prior}
In this subsection, we propose applying the pushforward measure from the f-CNF model as a suitable prior over infinite-dimensional spaces ($\mu_{\textbf{prior}}=\mu_{f_{\theta}}$, where $\mu_{f_{\theta}}$ is the transformed measure illustrated in Subsection \ref{subsec:Basic setting}). By training the f-CNF on the existing historical dataset, we can encode the historical information into the prior (we refer to it as the f-CNF prior in the remainder of this article). In this configuration, the Bayesian prior effectively encodes historical information, while the likelihood component $\mathbb{L}(\bm{d}|u)$ of the Bayes' formula accounts for the measurement data from the inverse problem. Leveraging Theorem $\ref{thm:RN derivative}$, we impose the RN estimation of the f-CNF prior:
$
\frac{d\mu_{f_{\theta}}}{d\mu_0}(u) = \exp(\mathcal{W}_{\theta}(u)).
$
Under this choice, Bayes’ formula becomes
\begin{align}\label{equ:bayes_f_CNF_paper}
\frac{d\mu}{d\mu_{f_{\theta}}}(u) = \frac{1}{Z_{\mu}^{\theta}} \exp(-\Phi(\bm{d}|u)).
\end{align}
where $Z_{\mu}^{\theta}$ is a positive finite constant given by $Z_{\mu}^{\theta} = \int_{\mathcal{H}_u} \exp (-\Phi(\bm{d}|u))\mu_{f_{\theta}}(du)$.

Recovering model parameters requires posterior inference from \eqref{equ:bayes_f_CNF_paper}, which is routinely performed using Markov chain Monte Carlo (MCMC) type algorithms. However, since classical infinite-dimensional MCMC sampling algorithm (e.g., pCN \cite{cotter2013,dashti2013bayesian}) are typically based on a Gaussian prior, it is necessary to express the Bayesian formulation \eqref{equ:bayes_f_CNF_paper} in terms of the pre-transformed Gaussian measure $\mu_0$ of the f-CNF as

\begin{align}\label{equ:bayesnewprior_paper}
\frac{d\mu}{d\mu_{0}}(u) = \frac{1}{Z_{\mu}^{\theta}} \exp\left(-\Phi(\bm{d}|u)+\mathcal{W}_{\theta}(u)\right).
\end{align}

It can be observed from formula \eqref{equ:bayesnewprior_paper} that, compared with directly applying measure $\mu_0$ as the prior, formulation \eqref{equ:bayesnewprior_paper} yields a new prior distribution transformed from $\mu_0$ via the trained f-CNF model $f_\theta$, where the transformation is implicitly reflected by the term $\mathcal{W}_\theta(u)$. This term directly depends on $f_\theta$ and can effectively incorporate useful information from the historical dataset. 

Our f-CNF prior exhibits a framework similar to the TV-Gaussian prior \cite{tvgaussian}. The TV-Gaussian prior constructs its RN density using a handcrafted TV regularizer
$
\frac{d\mu_{\rm pr}}{d\mu_0}(u) \propto \exp\bigl(-\lambda\|u\|_{\rm TV}\bigr),
$
while our prior adopts $\mathcal{W}_{\theta}(u)$ coupled with the neural flow $f_\theta$ for density characterization. Specifically, the TV-Gaussian prior focuses on accurately fitting functions with local jumps and sharp discontinuities, while our data-driven f-CNF prior is able to learn general prior knowledge from diverse historical dataset. Furthermore, when applying the pCN algorithm, our f-CNF-based posterior can be accelerated using the splitting strategy widely adopted in the TV-Gaussian framework \cite{tvgaussian}, see Appendix \ref{subappendix:Splitting acceleration strategy} for details.
 

Although the pCN algorithm provides accurate posterior approximations with sufficient samples \cite{pinski2015algorithms,Pinski2015SIAMMA,cotter2013} and is directly applicable via Bayes' formula \eqref{equ:bayesnewprior_paper}, it suffers from long mixing times, particularly with expensive forward solvers. Therefore, motivated by \cite{lu2025sequential}, we adopt a tailored sequential Monte Carlo with Gaussian mixture (SMC-GM) algorithm for our f-CNF prior to accelerate sampling.
 
 SMC-GM is an efficient sequential Monte Carlo algorithm for infinite-dimensional Bayesian inverse problems. It adopts Gaussian mixture mutations to replace standard MCMC updates, adjusts particle weights via potential functions, and conducts resampling to avoid particle degeneracy. With a tempering schedule that gradually incorporates measurement information, the algorithm iteratively transforms prior samples from $\mu_{f_{\theta}}$ into posterior samples. We refer readers to \cite{lu2025sequential} for complete algorithmic extensions and theoretical analysis, and provide full implementation details of the pCN and SMC-GM algorithms in sections \ref{subappendix:pCN Sampling under f-CNF Priors}
and \ref{subappendix:SMC-GM algorithm} of the Appendix.
 
\subsection{Training with Direct Dataset}\label{subsec:direct data}
The previous section discusses posterior sampling with the f-CNF prior. This and subsequent subsections will present two f-CNF prior learning strategies tailored to different data scenarios.  Specifically, this subsection focuses on the learning approach using dataset of the model parameters.

For example, consider a large historical dataset $\mathbb{D}_1 = \{u_n\}_{n=1}^{N_{\textbf{total}}}$ defined on $\mathcal{H}_u$, where each $u_n$ is a sample drawn from the data-generating measure $\mu_{\textbf{data}}$ of the model parameters. To construct a data-informed prior from these empirical samples, we adopt a generative flow model. Concretely, we train the f-CNF pushforward measure $\mu_{f_\theta}$ such that $\mathbb{D}_1$ approximately follows $\mu_{f_\theta}$. This embeds all information from $\mathbb{D}_1$ into the learned prior $\mu_{f_\theta}$, which can later serve as a valid Bayesian prior for newly observed data $\bm{d}$. The training procedure mirrors the standard setting of continuous normalizing flows for generative modeling on Euclidean spaces \cite{NODE}.
 
We define the direct-data training objective using the pointwise representative of the log-RN density:
\[
\mathcal{L}_1(\theta)
=-\mathbb{E}_{\mu_{\textbf{data}}}
\left[\mathcal{W}_{\theta}(u)\right]
\approx \widehat{\mathcal{L}}_1(\theta)
=-\frac{1}{N_{\textbf{total}}}
\sum_{n=1}^{N_{\textbf{total}}}\mathcal{W}_{\theta}(u_n).
\]

The specific training algorithm can be found in Algorithm \ref{alg:generate model}.
\begin{algorithm}[ht!]
\caption{f-CNF for Direct Data}
\label{alg:generate model}
\begin{algorithmic}[1]
\STATE Initialize parameters $\theta \leftarrow \theta_0$, batch size $N$, iteration number $K$, learning rate schedule $\alpha_k$, and set iteration counter $k \leftarrow 0$.
\REPEAT
    \STATE Sample a mini-batch of $N$ functions $\{u_i\}_{i=1}^{N}$ 
           from the full dataset $\mathbb{D}_1 = \{u_n\}_{n=1}^{N_{\text{total}}}$.
    \STATE Update parameters via gradient descent:
           \[
               \theta_{k+1} = \theta_k - \alpha_k \nabla_{\theta_k} \widehat{\mathcal{L}}_1(\theta_k).
           \]
           Stochastic gradient-based optimizers, such as Adam \cite{Adam}, 
           can be readily adopted in this step.
    \STATE Set $k \leftarrow k + 1$.
\UNTIL{$k = K$}
\STATE Return the learned parameter $\theta = \theta_K$.
\end{algorithmic}
\end{algorithm}

\subsection{Training with Indirect Dataset}\label{subsec:indirect data}
In prior distribution learning, true model parameters are generally unobservable in practical applications. In most cases, available historical records are physical measurements from real world processes, rather than direct samples of the unknown parameters of interest. To address this issue, various indirect prior learning methods have been developed to infer a prior from observational measurements instead of explicit parameter samples. In this subsection, we adopt the framework and data construction strategy in \cite{akyildiz2025efficient} to develop a practical application of our f-CNF model.

The work in \cite{jia2026nonparametric, akyildiz2025efficient} considers settings with abundant historical measurement data aggregated as $\mathbb{D}_2=\{\bm{d}_n\}_{n=1}^{N_{\text{total}}}$. All measurements in $\mathbb{D}_2$ are generated via the observation model $\bm{d}=\mathcal{S}_{\textbf{data}}\mathcal{G}(u)+\bm{\epsilon}_{\textbf{data}}$, where $\bm{\epsilon}_{\textbf{data}}$ denotes the observation noise and $\mathcal{S}_{\textbf{data}}$ represents the observation operator acting on the solution of the PDE. Specifically, there exists a latent set of model parameters $\mathbb{D}_{\textbf{hid}}=\{u_i\}_{i=1}^{N_{\text{total}}}$ such that:
$\mathbb{D}_2=\{\bm{d}_i \mid \bm{d}_i=\mathcal{S}_{\textbf{data}}\mathcal{G}(u_i)+\bm{\epsilon}_{\textbf{data}},\, u_i \in \mathbb{D}_{\textbf{hid}}\}.$ Evidently, all accessible information about $\mathbb{D}_{\textbf{hid}}$ is contained within the observable dataset $\mathbb{D}_2$.
 Our objective is to train an f-CNF prior $\mu_{f_{\theta}}$ using the observational information from $\mathbb{D}_2$, such that $\mu_{f_{\theta}}$ can accurately approximate the latent dataset $\mathbb{D}_{\textbf{hid}}$.

In such scenarios, the f-CNF prior $\mu_{f_{\theta}}$ can no longer be trained using the negative log-likelihood loss. Here, similar to \cite{akyildiz2025efficient}, our loss builds on the Sinkhorn distance, an entropy-regularized Wasserstein variant devised to bypass the original Wasserstein’s computational bottlenecks and provide efficient differentiable loss for deep learning. The specific form of the loss function is as follows:
$$
\mathcal{L}_2(P, Q) = W_\alpha^2(P, Q) = \min_{\pi \in \Pi(P, Q)} \mathbb{E}_{(\boldsymbol{x}, \boldsymbol{y}) \sim \pi} \left[ \|\boldsymbol{x} - \boldsymbol{y}\|^2_{\mathbb{R}^{N_d}} \right] + \alpha H(\pi).
$$
Here, \(P\) and \(Q\) represent two distributions for which only samples are available, while their exact densities remain unknown. $H(\pi)$ denotes the entropy of the transport plan $\pi$. For further details on this distance, please refer to \cite{feydy2019interpolating}.

Under assumptions of this subsection, historical parameter information is encoded via abundant measurements $\mathbb{D}_2=\{\bm{d}_n\}_{n=1}^{N_{\text{total}}}$.
 Let $\nu$ denote the law of $\{\bm{d}_n\}_{n=1}^{N_{\text{total}}}$. Suppose we have $\bar{\bm{d}}_n = \mathcal{S}_{\textbf{data}}\mathcal{G}(u_n) + \bm{\epsilon}_{\textbf{data}}$, where $\{u_n\}_{n=1}^{N_{\text{total}}}$ are samples from $\mu_{f_{\theta}}$, and let $\bar{\nu}$ denote the law of $\{\bar{\bm{d}}_n\}_{n=1}^{N_{\text{total}}}$. Our objective is to ensure that $\bar{\nu}$ is close to $\nu$ in some appropriate sense. It is worth noting that while the explicit forms of measures $\nu$ and $\bar{\nu}$ remain unknown, their associated samples are readily available.
 In practice, such measures are often represented by their empirical counterparts. To effectively and efficiently compare these empirical measures, we employ $\mathcal{L}_2$ as the loss function for training.

Let $*$ denote the convolution of measures. To quantify the discrepancy between probability measures in the measurement space, we introduce a distance $\mathcal{L}_2$ and define the corresponding loss function in this subsection as follows:
\[
\mathcal{J}(\theta) =\mathcal{L}_2( \nu, \bar{\nu} )= \mathcal{L}_2\Bigl( \nu, \eta * (\mathcal{S}_{\textbf{data}}\mathcal{G})_\#\mu_{f_{\theta}} \Bigr),
\]
where $\eta$ denotes the Gaussian distribution of the data noise $\bm{\epsilon}_{\textbf{data}}$, i.e., $\bm{\epsilon}_{\textbf{data}} \sim \eta$. Detailed theoretical analysis of the loss function can be found in \cite{akyildiz2025efficient}, and the specific training algorithm is presented in Algorithm \ref{alg:generate model indirect}.
\begin{algorithm}[ht!]
\caption{f-CNF for Indirect Data}
\label{alg:generate model indirect}
\begin{algorithmic}[1]
\STATE Initialize parameters $\theta \leftarrow \theta_0$, batch size $N$, 
       total iterations $K$, learning rate schedule $\{\alpha_k\}$, 
       and set iteration counter $k \leftarrow 0$.
\REPEAT
    \STATE Draw $N$ functions $\{u_i\}_{i=1}^{N}$ from $\mu_{f_{\theta_k}}$.
    \STATE For $i = 1, \ldots, N$, compute the  observations 
           $\bar{\bm{d}}_i = \mathcal{S}_{\textbf{data}}\mathcal{G}(u_i) + \bm{\epsilon}_{\textbf{data}}$.
    \STATE Sample a mini batch of $N$ measurement  $\{\bm{d}_i\}_{i=1}^{N}$ from the  dataset $\{\bm{d}_i\}_{i=1}^{N_{\text{total}}}$.
    \STATE Compute the loss $\mathcal{J}(\theta_k) = \mathcal{L}_2(\nu, \bar{\nu})$ by approximating the exact expectations with the empirical means computed from $\{\bm{d}_i\}_{i=1}^{N}$ and $\{\bar{\bm{d}}_i\}_{i=1}^{N}$.

    \STATE Update parameters via gradient descent:
$$
               \theta_{k+1} = \theta_k - \alpha_k \nabla_{\theta_k} \mathcal{J}(\theta_k).
$$
           Stochastic gradient-based optimizers, such as Adam \cite{Adam}, 
           can be readily adopted in this step.
    \STATE Set $k \leftarrow k + 1$.
\UNTIL{$k = K$}
\STATE \textbf{Return} the learned parameters $\theta_K$.
\end{algorithmic}
\end{algorithm}

It is noteworthy that, compared to training the f-CNF prior using extensive historical dataset of model parameter as discussed in Subsection \ref{subsec:direct data}, learning the prior from historical measurement data $\{\bm{d}_n\}_{n=1}^{N_{\text{total}}}$ is often considerably more challenging. Specifically, consider the forward problem $\bm{d} = \mathcal{S}_{\textbf{data}}\mathcal{G}(u) + \bm{\epsilon}_{\textbf{data}}$. When excessive information loss occurs in the forward PDE process (associated with $\mathcal{G}$), the number of measurement points is insufficient (associated with $\mathcal{S}_{\textbf{data}}$), or the measurement noise is excessively large (associated with $\bm{\epsilon}_{\textbf{data}}$), the information about $u$ contained in the dataset $\{\bm{d}_n\}_{n=1}^{N_{\text{total}}}$ diminishes, thereby introducing significant difficulties in training the parameters of the f-CNF. 

\subsection{Other Settings}\label{subsec:other setting}
In the previous two subsections, we present the applications of the proposed f-CNF model with either a large set of model parameters $\mathbb{D}_1=\{u_n\}_{n=1}^{N_{\text{total}}}$ or abundant measurement data $\mathbb{D}_2=\{\bm{d}_n\}_{n=1}^{N_{\text{total}}}$. In fact, the f-CNF prior is applicable to broader scenarios and can be combined with various existing methods. For example, it can be embedded into meta-learning paradigms. As studied in \cite{jia2026nonparametric}, an infinite-dimensional PAC-Bayesian framework is established to learn data-adaptive priors, which yields more task-specified priors and improves sampling efficiency. Nevertheless, their work only focuses on Gaussian priors. Substituting the Gaussian prior with our f-CNF prior in such a framework constitutes a valuable research direction for future investigation. 

In addition, the f-CNF model can be associated with functional flow matching
\cite{kerrigan2023functional, shikhman2026discretization}, which extends
flow matching to infinite-dimensional Hilbert spaces without relying on Lebesgue densities. Once f-CNF is trained in this way, by Theorem~\ref{thm:RN derivative}, the measures before and after the learned Neural ODE flow are equivalent, and the corresponding RN derivative is available for Bayesian analysis.

Our model also exhibits connections to diffusion models, which we briefly discuss in what follows.
Infinite-dimensional score based diffusion models are well-established on infinite-dimensional spaces \cite{diffusionf2}. However, these models are formulated through stochastic differential
equations (SDEs), which do not directly provide explicit densities
relative to a reference measure. This complicates their use in
Bayesian inference algorithms requiring density evaluation, as
in their finite-dimensional counterparts
\cite{song2021maximum, feng2023score}. A research for finite-dimensional diffusion models \cite{feng2023score} shows that the SDE of a diffusion model can be reduced to a deterministic probability flow ODE. Interpreting this ODE as a continuous normalizing flow allows one to apply the change-of-variable theorem to compute the exact log-density, which enables diffusion models to act as rigorous Bayesian priors. This approach provides a guideline for developing RN estimation for diffusion models on infinite-dimensional spaces.

In the infinite-dimensional setting \cite{diffusionf2}, the generative model is governed by a reverse SDE that maps a Gaussian distribution back to the data distribution:
\begin{align}\label{equ:diffusion_infinite}
du_t = \frac{1}{2} u_t dt + s_{\theta}(t, u_t) dt + dW_t^{\mathcal{C}},
\end{align}
where $s_{\theta}(t,u)$ denotes the approximate score function, $W_t^{\mathcal{C}}$ denote the $\mathcal{C}$-Wiener processes, and $u_0 \sim \mathcal{N}(0,\mathcal{C})$ is the initial Gaussian measure. Analogous to \cite{feng2023score}, we hypothesize that the Brownian motion term $\mathrm{d}W_t^{\mathcal{C}}$ can be eliminated while preserving marginal distributions. This yields the probability flow ODE associated with the infinite-dimensional forward process:
\begin{align}\label{equ:diffusion_ODE_infinite}
\frac{\mathrm{d}u_t}{\mathrm{d}t} = f_{\theta}(u_t,t).
\end{align}
The conjecture is established strictly following the framework in \cite{feng2023score}, and the explicit form of $f_{\theta}(u_t,t)$ for the reverse SDE \eqref{equ:diffusion_infinite} requires further investigation. Such deterministic flow defines an f-CNF model. Even though Lebesgue densities are not well-defined in infinite-dimensional spaces, this model admits a rigorous RN derivative as the generalized density relative to the Gaussian reference measure $\mathcal{N}(0,\mathcal{C})$, whenever the ODE \eqref{equ:diffusion_ODE_infinite} satisfies the assumptions stated in Theorem \ref{thm:RN derivative}.

In summary, the f-CNF model proposed here is not limited to the two prior construction methods presented in this paper. Its integration with prior learning approaches in broader contexts remains a subject for future research.

\section{Well-Posedness of Bayesian Inverse Problems}\label{sec:wellposedness}
In the Bayesian approach to inverse problems, the uncertain model parameter is represented as a random variable governed by a prior measure that encodes initial uncertainty. The observation of data constitutes an event upon which this prior is conditioned. The solution to the Bayesian inverse problem is the resulting posterior measure, defined as this conditional probability distribution. A key concern is whether small perturbations in observational data induce significant changes in the posterior measure, namely, whether the posterior is stable with respect to data perturbations. Reference \cite{MR4624340} has investigated this issue under various settings. Based on the results established in \cite{MR4624340}, we further analyze the well-posedness of the Bayesian inverse problem equipped with the newly proposed f-CNF prior. Consistent with the framework in \cite{MR4624340}, we introduce the following definitions.
\begin{definition}\label{def:wellposedness}Assume that all the probability measures on $\mathcal{H}_u$ combines the set $\mathbb{P}$, and $\rho$ is a metric on $\mathbb{P}$. The Bayesian inverse problem with respect to \eqref{equ:bayes infinite} is $\rho$-well-posed if
\begin{itemize}
    \item Bayes' formula \eqref{equ:bayes infinite} has a unique posterior in $\mathbb{P}$ (existence and uniqueness).
    \item The posterior of \eqref{equ:bayes infinite} depends continuously on the data $\bm{d}$, i.e., $(Y, \|  \cdot  \|_{Y}) \ni \bm{d} \rightarrow  \mu_{\bm{d}} \in  (\mathbb{P}, \rho)$ is a continuous function (stability).
\end{itemize}
\end{definition}

It is noteworthy that the definition of well-posedness in Definition \ref{def:wellposedness} concerns the continuity of the measures with respect to the observed data $\bm{d}$. Consequently, by employing different metrics $\rho$ between probability measures, we can formulate distinct notions of well-posedness. Following the exposition in \cite{MR4624340}, we present several classical definitions of well-posedness below.
\begin{definition}
We classify $(\mathbb{P},\rho)$ well-posedness according to the choice of the metric $\rho$ as follows:
\begin{itemize}
\item weak well-posedness, if $\rho$ is the Prokhorov metric;
\item total variation well-posedness, if $\rho$ is the total variation distance;
\item Hellinger well-posedness, if $\rho$ is the Hellinger distance;
\item Wasserstein(p) well-posedness, if $\rho$ is the Wasserstein(p) distance.
\end{itemize}
\end{definition}

For the Bayesian formulation associated with the f-CNF prior $\mu_{f_{\theta}}$ constructed herein, as given in formula \eqref{equ:bayes_f_CNF_paper}, Bayesian well-posedness can be established provided that the likelihood function satisfies certain suitable conditions.
\begin{theorem}\label{thm:wellposedness}
    Consider the Bayesian problem with f-CNF prior with the following assumptions hold:
    \begin{itemize}
        \item the ODE corresponding to f-CNF satisfised the conditions of Theorem \ref{thm:RN derivative};
        \item The likelihood component $\mathbb{L}(\cdot|u)$ is continuous;
        \item There exist a constant $M$ independent of $u$ and $\bm{d}$ with  $0<\mathbb{L}(\bm{d}|u) \leq M$;
    \end{itemize}
    then the Bayesian inverse problem is weakly, Hellinger, total variation, and Wasserstein(p) well-posed.
\end{theorem}

The detailed proof is provided in the Appendix \ref{subappendix:proof of wellposedness}.

\section{Numerical Experiments}\label{sec:numerical experiment}

In the study \cite{NODE}, finite-dimensional continuous planar flows were introduced as canonical continuous normalizing flows; this construction extends naturally to the infinite-dimensional setting of the present work:
\begin{align}\label{equ:experiment}
\left\{
\begin{array}{l}
\displaystyle \frac{dv(t)}{dt} = \sum_{i=1}^{L} u_i(t)\sigma(\langle v(t),w_i(t)\rangle_{\mathcal{H}_u} + b_i(t)), \\
v(0) = u_0,
\end{array}
\right.
\end{align}
where, for any fixed $t\in[0,T]$ and $i=1,\ldots,L$, $u_i(t),w_i(t) \in \mathcal{H}_u$, $b_i(t) \in \mathbb{R}$, and $\sigma(x)=\tanh(x)$. Parameterizing the time-dependent functions $u_i(t)$, $w_i(t)$, and $b_i(t)$ is essential for computational efficiency and satisfying Theorem \ref{thm:RN derivative}. Let $\{\phi_j\}_{j=1}^{M}$ be the first $M$ eigenfunctions of the covariance operator associated with $\mu_0$ in Theorem \ref{thm:RN derivative}, where $M$ is the internal basis truncation dimension that fixes the model capacity and remains constant after training. We expand $u_i(t)$ and $w_i(t)$ as
\begin{align}\label{equ:un}
    u_i(t)=\sum\limits_{j=1}^{M}\alpha_i^{(j)}(t) \phi_j, \qquad
    w_i(t)=\sum\limits_{j=1}^{M}\beta_i^{(j)}(t) \phi_j.
\end{align}

For brevity, denote $\bm{\alpha}_i(t)=(\alpha_i^{(1)}(t),\dots,\alpha_i^{(M)}(t))$ and $\bm{\beta}_i(t)=(\beta_i^{(1)}(t),\dots,\beta_i^{(M)}(t))$. Since $\bm{\alpha}_i, \bm{\beta}_i: \mathbb{R} \to \mathbb{R}^M$ and $b_i: \mathbb{R} \to \mathbb{R}$ share the real-valued input $t \in \mathbb{R}$ and their combined output lies in $\mathbb{R}^{2ML+L}$, these mappings can be jointly parameterized by a fully connected neural network with trainable parameters $\theta$.

A key feature of this eigenbasis construction is its mesh-independent parameterization: since each basis function $\phi_j$ can be evaluated pointwise at different spatial locations, a trained f-CNF with fixed $M$ can be applied to input functions represented on different spatial meshes without retraining the network parameters. This is in the same spirit as Fourier neural operators (FNO) \cite{li2020fourier}, where a fixed number of internal spectral modes enables mesh-invariant evaluation. We emphasize that, similar with FNO, this property concerns only the discretization of the function domain. It does not imply invariance to the internal truncation dimension $M$. As with any spectral method, $M$ controls model expressivity and should be chosen appropriately for the complexity of the target problem.

We verify that the network \eqref{equ:experiment} constructed via the aforementioned approach satisfies the conditions of Theorem \ref{thm:RN derivative} (see Appendix \ref{subappendix:Verification of the Continuous Planar Flow Conditions} for details). Other ODE-driven f-CNF architectures satisfying Theorem \ref{thm:RN derivative} are given in Appendix \ref{subappendix:Other flow model}. Throughout the numerical experiments, we adopt the formulation in \eqref{equ:experiment}. Unless noted otherwise, model \eqref{equ:un} will employ a fully connected neural network with three hidden layers (12 neurons per layer), $\tanh(\cdot)$ activation function, and fixed $M=20$, $L=12$.


The code used for simulations is available on GitHub\footnote{\url{https://github.com/jjx323/FunctionalContinuousNormalizingFlow}}. All programs ran on a system with a 13th Gen Intel(R) Core(TM) i7-13700, NVIDIA GeForce RTX 4090, and Ubuntu 20.04.5 LTS.

\subsection{Simple Smooth}\label{subsec:simplesmooth}
Consider the inverse source problem studied in \cite{sui2024noncentered}, governed by the following elliptic equation:
\begin{align}\label{equ:simple}
\begin{split}
 -\alpha \Delta w + w &= u, \quad \text{in}\ \Omega, \\ 
 \frac{\partial w}{\partial \bm{n}} &= 0, \quad \text{on}\ \partial \Omega,
\end{split}
\end{align}
where $\Omega = (0, 1) \subset \mathbb{R}$, $\alpha > 0$ is a positive constant, and $\bm{n}$ denotes the outward unit normal vector. The forward operator is defined by
\begin{equation*}
\mathcal{S}\mathcal{G}u =(w(x^1), w(x^2), \ldots, w(x^{N_{\bm{d}}}) )^T,
\end{equation*}
where $\mathcal{G}$ denotes the PDE solution operator from $\mathcal{H}_u:=L^2(\Omega)$ to $\mathcal{H}_{w}:=H^2(\Omega)$, $\mathcal{S}$ denotes the measurement operator from $\mathcal{H}_{w}$ to $\mathbb{R}^{N_{\bm{d}}}$, $u \in \mathcal{H}_u$ is the source function, $w\in \mathcal{H}_w$ is the solution to the elliptic equation \eqref{equ:simple}, and $x^i \in \Omega$ for $i = 1, 2, \ldots ,  N_{\bm{d}}$.
The data $\bm{d}$ are related to $u$ through the forward model
\begin{align}\label{equ:simple inverse}
 \bm{d} = \mathcal{S}\mathcal{G}u + \bm{\epsilon}.
\end{align}
In this subsection, the noise $\bm{\epsilon}$ is taken to be 10$\%$ Gaussian random noise, i.e., $\bm{\epsilon} \sim \mathcal{N}(0, \bm{\Gamma}_{\text{noise}})$, where $\bm{\Gamma}_{\text{noise}} = \tau^{-1}\textbf{I}$ and $\tau^{-1} = (0.1 \lVert \mathcal{S}\mathcal{G}u\rVert_{\infty})^2$.

By applying the Bayesian framework introduced in Subsection \ref{subsec:Bayesian Inverse Problem}, we reformulate the inverse problem under consideration as a posterior inference problem on infinite-dimensional space. In this setting, an appropriate prior must be specified in the Bayes' formula \eqref{equ:bayes infinite}. A conventional choice is to employ a weakly informed Gaussian prior $\mu_0 = \mathcal{N}(0,\mathcal{C}_0)$, as considered in \cite{dashti2013bayesian,stuart2010inverse,cotter2013}. However, such priors may fail to produce satisfactory inversion results when the measurement information contained in the likelihood function of \eqref{equ:bayes infinite} is limited. To address this issue, we employ the f-CNF model presented in Sections \ref{sec:f-CNF} and \ref{sec:IP_and_Prior}, where the weakly informed reference prior $\mu_0$ is pushed forward through a trained f-CNF network $f_{\theta}$, yielding a problem-adapted prior $\mu_{f_{\theta}}$ that is better suited to the target inverse problem. Here, the covariance operator of the baseline Gaussian measure $\mu_0$ is defined as $\mathcal{C}_0 = \big(\mathrm{I} - \alpha \Delta\big)^{-2}$, where $\alpha = 0.05$ and $\Delta$ denotes the Laplace operator on $\Omega$ with homogeneous Neumann boundary conditions.

To train the f-CNF network, we construct direct dataset  $\mathbb{D}_{\mathbf{train}}^1$ and indirect dataset $\mathbb{D}_{\mathbf{train}}^2$ similar to those in \cite{li2020fourier, akyildiz2025efficient}, with the detailed definitions given below:
\begin{itemize}
    \item Let $u = \cos(2\pi x)+\cos(4\pi x)+0.3\sum\limits_{k=1}^{20} \xi_k \cos(k\pi x) /(1+k^2) $, where $\xi_k \sim \mathcal{N}(0,1)$. The direct training dataset of model parameters is then defined as $\mathbb{D}_{\textbf{train}}^1=\{u_i\}_{i=1}^{5000}$, where the samples are generated according to the above sampling strategy.
    \item Define the indirect training dataset by $\mathbb{D}_{\textbf{train}}^2=\{\bm{d}_i \mid \bm{d}_i=\mathcal{S}_{\textbf{data}}\mathcal{G}(u_i)+\bm{\epsilon}_{\textbf{data}},\, u_i \in \mathbb{D}_{\textbf{train}}^1\}$, where the measurement points of $\mathcal{S}_{\textbf{data}}$ are specified as $\lbrace w(x^i) \mid i = 1, 2, \ldots, 10 \rbrace$ with $x^{i}=i/10$, and $\bm{\epsilon}_{\textbf{data}}$ denotes Gaussian white noise with variance $0.01$.
\end{itemize}

\begin{figure}[ht!]
	\centering
	
	\subfloat[Samples ($\mu_0$)]{
		\includegraphics[ keepaspectratio=true, width=0.32\textwidth, clip=true]{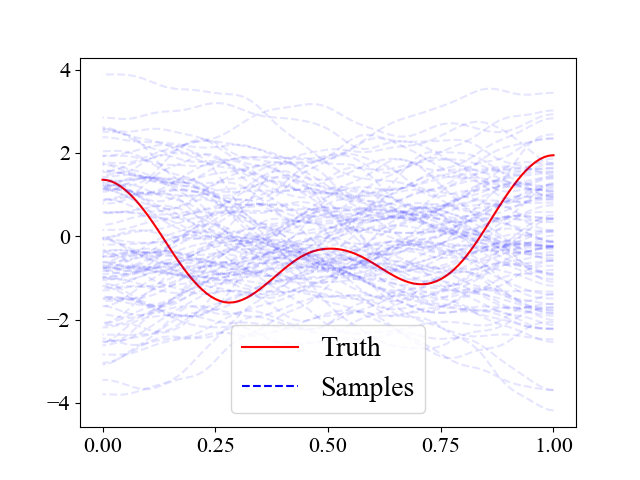}} 
	\subfloat[Samples ($\mu_{f_{\theta}^1}$)]{
		\includegraphics[ keepaspectratio=true, width=0.32\textwidth, clip=true]{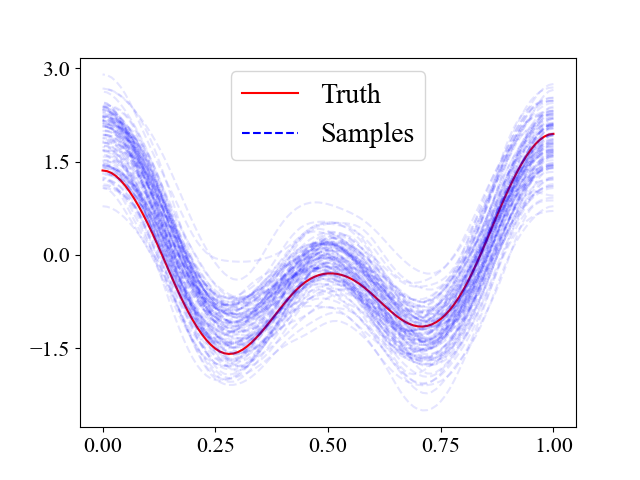}} 
        \subfloat[Samples ($\mu_{f_{\theta}^2}$)]{
		\includegraphics[ keepaspectratio=true, width=0.32\textwidth, clip=true]{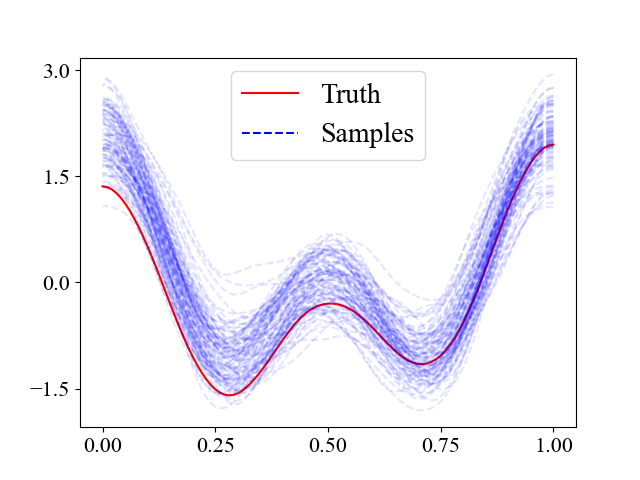}} 
	\caption{\emph{\small  Samples from $\mu_0$ and trained f-CNF priors $\mu_{f_{\theta}^1}$ and $\mu_{f_{\theta}^2}$. The red solid line represent the background truth $u^{\dagger}$. (a): Samples from $\mu_0$. (b): Samples from $\mu_{f_{\theta}^1}$ trained with direct dataset. (c):Samples from $\mu_{f_{\theta}^2}$ trained with indirect dataset. }}
       \label{pic:simplesmooth prior}
\end{figure}

Subsequently, the ground truth $u^{\dagger}$ of the model \eqref{equ:simple inverse} is generated using the same sampling strategy as that used for the training dataset $\mathbb{D}_{\textbf{train}}^1$. To better illustrate the advantages of the learned f-CNF prior $\mu_{f_{\theta}}$ over $\mu_0$, the measurement points associated with the measurement operator $\mathcal{S}$ in \eqref{equ:simple inverse} are chosen to be relatively sparse, namely, $\{x^i\}_{i=1}^{5}$ with $x^i = 0.1 + 0.2(i-1)$. The corresponding measurement data are then generated by
$$
\bm{d}^{\dagger}=\mathcal{S}\mathcal{G}(u^{\dagger})+\bm{\epsilon}.
$$
We train the f-CNF prior $\mu_{f_{\theta}^1}$ on the direct dataset $\mathbb{D}_{\mathrm{train}}^1$ via Algorithm \ref{alg:generate model}. Likewise, $\mu_{f_{\theta}^2}$ is trained on the indirect dataset $\mathbb{D}_{\mathrm{train}}^2$ using Algorithm \ref{alg:generate model indirect}. Samples are then drawn from these two learned priors and the baseline measure $\mu_0$ for comparison against the ground truth $u^\dagger$. The corresponding results are presented in Figure \ref{pic:simplesmooth prior}. It can be observed that samples drawn from $\mu_{f_{\theta}^1}$ and $\mu_{f_{\theta}^2}$ exhibit greater similarity to $u^{\dagger}$ than those generated from $\mu_0$.

Subsequently, we adopt $\mu_0$, $\mu_{f_{\theta}^1}$, and $\mu_{f_{\theta}^2}$ as the prior measures in the Bayes' formula \eqref{equ:bayes infinite}. For this simple inverse problem, we generate posterior samples via the pCN algorithm. The step size is set to $\beta=0.05$ for $\mu_0$, and $\beta=0.02$ for $\mu_{f_{\theta}^1}$ and $\mu_{f_{\theta}^2}$ (note that the prior $\mu_0$ admits a larger step size than $\mu_{f_{\theta}}$; the underlying reasons are elaborated in the Appendix \ref{subappendix:Splitting acceleration strategy}). For all three priors, we run $5 \times 10^7$ iterations and discard the first $1 \times 10^6$ samples as burn-in.

\begin{figure}[ht!]
	\centering
	
	\subfloat[Posterior ($\mu_0$)]{
		\includegraphics[ keepaspectratio=true, width=0.32\textwidth, clip=true]{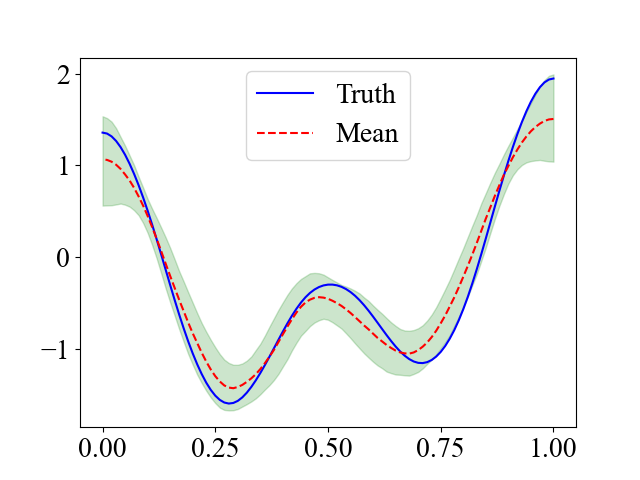}} 
	\subfloat[Posterior ($\mu_{f_{\theta}^1}$)]{
		\includegraphics[ keepaspectratio=true, width=0.32\textwidth, clip=true]{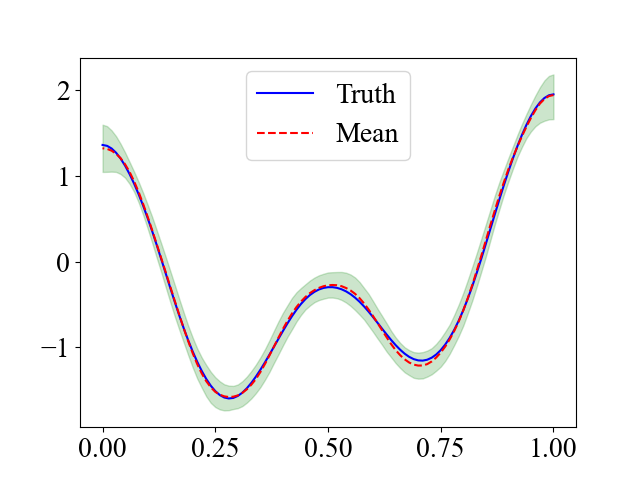}} 
        \subfloat[Posterior ($\mu_{f_{\theta}^2}$)]{
		\includegraphics[ keepaspectratio=true, width=0.32\textwidth, clip=true]{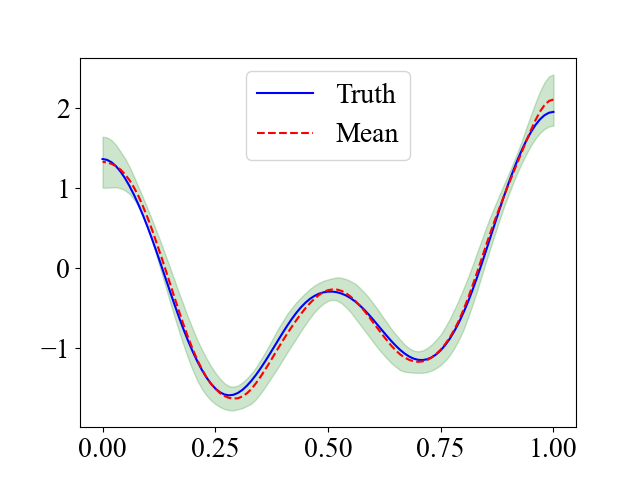}} 
	\\
	\subfloat[Covariance ($\mu_0$)]{
		\includegraphics[ keepaspectratio=true, width=0.32\textwidth, clip=true]{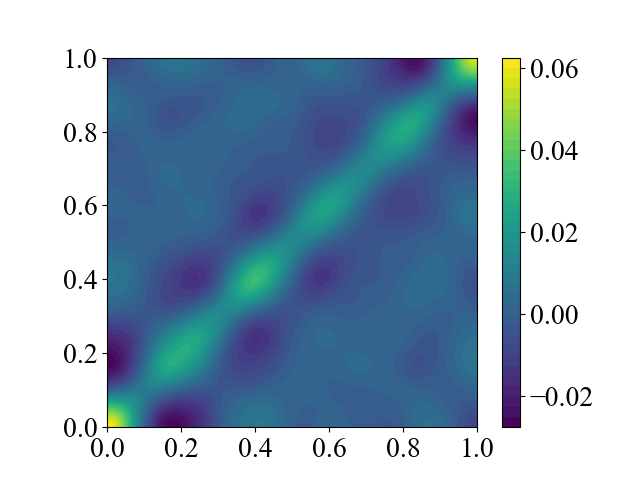}} 
	\subfloat[Covariance ($\mu_{f_{\theta}^1}$)]{
		\includegraphics[ keepaspectratio=true, width=0.32\textwidth, clip=true]{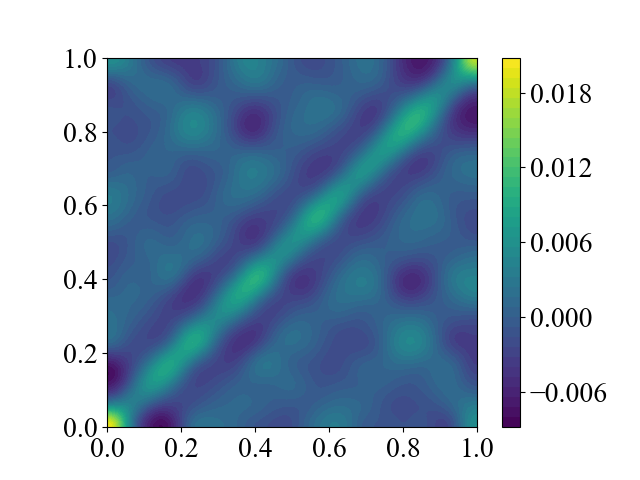}} 
        \subfloat[Covariance ($\mu_{f_{\theta}^2}$)]{
		\includegraphics[ keepaspectratio=true, width=0.32\textwidth, clip=true]{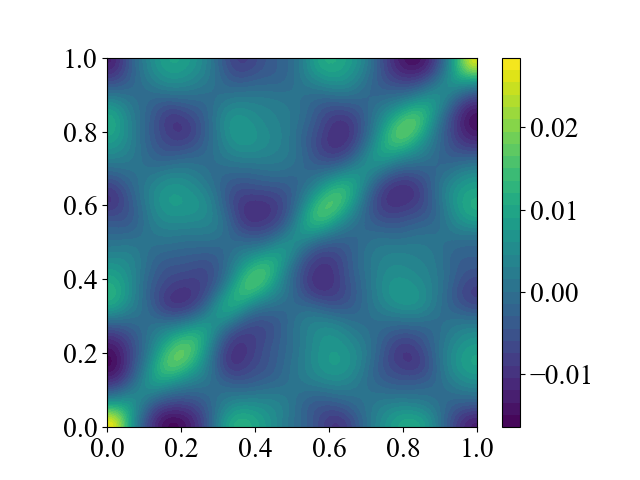}} 
	\caption{\emph{\small  The comparison of the Bayesian posteriors with $\mu_0$, $\mu_{f_{\theta}^1}$, $\mu_{f_{\theta}^2}$ priors. (a)(b)(c): The corresponding Bayesian posteriors with $\mu_0$, $\mu_{f_{\theta}^1}$, $\mu_{f_{\theta}^2}$ priors. The blue solid line denotes the background truth $u^{\dagger}$, the red dashed line represents the mean of the posterior, and the green shaded region indicates the $95\%$ credible interval of the posteriors.  (d)(e)(f): The covariance function of posteriors with $\mu_0$, $\mu_{f_{\theta}^1}$, $\mu_{f_{\theta}^2}$ priors.
 }}
\label{pic:simplesmooth posterior}
\end{figure}

For clarity, we recall the definition of the covariance function. Let $u$ be a random field defined on a domain $\Omega$, with mean function $\bar{u}$ and covariance function $c(x, y)$. Here $c(x, y)$ characterizes the covariance between $u(x)$ and $u(y)$:
\begin{align*}
	c(x, y) = \mathbb{E}((u(x)-\bar{u}(x))(u(y)-\bar{u}(y))), \quad \text{for} \ x, y \in \Omega.
\end{align*}
The corresponding covariance operator $\mathcal{C}$ is defined by
\begin{align*}
	(\mathcal{C}\phi^{\prime})(x) = \int_{\Omega}c(x, y)\phi^{\prime}(y)dy,
\end{align*}
where $\phi^{\prime}$ is a sufficiently regular function defined on $\Omega$. 
To visualize the Bayesian posterior, we plot its mean and credible region in Figures  \ref{pic:simplesmooth posterior}(a)–(c). The posterior covariance functions in panels (d)–(f) further quantify the uncertainty of the inversion results. Additionally, we compute the $L^2$ error between the posterior mean and the background truth, with results summarized in Table \ref{table:simplesmooth relative}.

\begin{table}[ht!]
	\renewcommand{\arraystretch}{1.5}
	\centering
		\caption{\emph{\small The relative errors between the background truth and the mean of posteriors.}} 
	\begin{tabular}{c|ccc}
		\hline $\text{Relative Error}$ & posterior ($\mu_0$)   & posterior ($\mu_{f_{\theta}^1}$) & posterior ($\mu_{f_{\theta}^2}$) \\
		\hline $\text{Background Truth}$ $u^{\dagger}$& $0.15641$& $0.04515$& $0.05693$ \\
		\hline
	\end{tabular}
 \label{table:simplesmooth relative}
\end{table}
From the numerical results in Table \ref{table:simplesmooth relative} and the graphical illustrations in Figure \ref{pic:simplesmooth posterior}, we observe that the posterior mean obtained by the Bayesian method with $\mu_0$ as the prior exhibits a relatively large error compared with the background truth $u^{\dagger}$. This indicates that the likelihood provides insufficient measurement information to accurately recover the true source function $u^{\dagger}$. In contrast, the proposed f-CNF priors $\mu_{f_{\theta}^1}$ and $\mu_{f_{\theta}^2}$, trained using direct and indirect data, respectively, mitigate this deficiency by incorporating additional prior information.

This phenomenon can be further understood by examining the covariance function of the Bayesian posterior. Subfigure (d) of Figure \ref{pic:simplesmooth posterior} shows that the covariance function of the posterior associated with the prior $\mu_0$ remains relatively large, indicating that sparse measurements lead to high uncertainty in the inversion results. In comparison, Subfigures (e) and (f) of Figure \ref{pic:simplesmooth posterior} show the posterior covariance functions corresponding to the Bayesian formulations equipped with the priors $\mu_{f_{\theta}^1}$ and $\mu_{f_{\theta}^2}$, respectively, where additional information is incorporated through the learned f-CNF priors. These learned priors enable the posterior distributions to produce narrower credible region and smaller covariance values. Moreover, since direct data inherently contain more information about the target parameter than indirect data, the prior $\mu_{f_{\theta}^1}$ yields better performance than $\mu_{f_{\theta}^2}$.

\subsection{Inverse Scattering Problem}\label{subsec:inverse_scar}

In this subsection, we investigate the application of the f-CNF prior to the inverse scattering problem \cite{scar1,scar2} using the available direct dataset $\mathbb{D}_{\textbf{train}}$. We first briefly introduce the sound-soft inverse scattering problem considered here. Let $\Omega\subset\mathbb{R}^2$ be a bounded, simply connected obstacle with boundary $\partial\Omega$. In $\mathbb{R}^2$, define $\mathbb{S}^1 = \{\xi \in \mathbb{R}^2 : |\xi| = 1\}$. Consider the illumination of $\Omega$ by the following plane wave:
\begin{equation*}
w^i(x) := e^{ikx \cdot \nu}, \quad x \in \mathbb{R}^2, \, \nu \in \mathbb{S}^1,
\end{equation*}
where $k > 0$ is the wavenumber and $\nu$ denotes the incident direction. Given the obstacle $\Omega$ and the incident direction $\nu$, the PDE solution operator determines the scattered field $w^s$, or equivalently the total field $w = w^i + w^s$, such that $w$ satisfies
\begin{align}
\Delta w + k^2 w &= 0, \quad \text{in } \mathbb{R}^2 \backslash \overline{\Omega}, \label{eq:helmholtz} \\
w &= 0, \quad \text{on } \partial\Omega, \label{eq:boundary} \\
\lim_{r \to \infty} \sqrt{r} &\left( \frac{\partial w^s}{\partial r} - ikw^s \right) = 0, \label{eq:sommerfeld}
\end{align}
where \eqref{eq:helmholtz} is the Helmholtz equation \cite{scar1,scar2,colton1996simple}, \eqref{eq:boundary} is the sound-soft boundary condition, and \eqref{eq:sommerfeld} is the Sommerfeld radiation condition.

Assume that the scatterer $\Omega$ is a star-shaped domain. The boundary $\Gamma$ of the scatterer $\Omega$ can be parameterized as
\begin{equation*}
\Gamma := r(\theta)(\cos\theta, \sin\theta) = \exp(u(\theta))(\cos\theta, \sin\theta), \quad \theta \in [0, 2\pi),
\end{equation*}
where $u(\theta) = \ln r(\theta)$ and $r(\theta)>0 $.
For a fixed wavenumber $k=5$, the incident field is defined as $w^i(x) := e^{ikx \cdot \nu}$. 

The forward problem is formulated by illuminating the obstacle with plane waves from eight distinct incident directions
\begin{align*}
\nu \in \left\{ \big(\cos\theta_i,\,\sin\theta_i\big) \,\bigg|\,
\theta_i = \frac{(i-1)\pi}{4},\ i=1,2,\dots,8 \right\}.
\end{align*}
For each incident direction, the corresponding scattered field is computed by solving the governing PDE that describes the interaction between the incident wave and the target obstacle.

For each incident direction, we collect phaseless noisy scattered-field measurements at three observation points $\{(6,0),(3\sqrt{3},3\sqrt{3}),(0,6)\}$. Let $\mathcal{G}_i$ denote the solution operator for the $i$-th incident wave, which maps $u \in \mathcal{H}_u$ to the scattered field $w_i \in \mathcal{H}_w$, i.e., $w_i = \mathcal{G}_i(u)$. Let $\mathcal{S}$ stands for the measurement operator that extracts phaseless data of the scattered-field at the three fixed points. The measurement model for the $i$-th incident wave is given by
$
\bm{d}_i = \mathcal{S}\mathcal{G}_i(u) + \bm{\epsilon}_i,
$
where $\bm{\epsilon}_i$ represents Gaussian noise.
By collecting the solution operators corresponding to all incident waves, we define the combined solution operator $\mathcal{G}$ by
$
\mathcal{G}(u)=(\mathcal{G}_1(u),\mathcal{G}_2(u), \ldots, \mathcal{G}_8(u)).
$  
The measurement vector can be written as
$$
\bm{d} = \mathcal{S}\mathcal{G}(u) + \bm{\epsilon},
$$
where $\bm{d}$ is a vector of length $3 \times 8$. The potential function in Bayes' formula \eqref{equ:bayes infinite} is defined as
\begin{equation*}
\Phi(u) := \frac{1}{2}\left\lVert \mathcal{S}\mathcal{G}(u) - \bm{d}\right\rVert_{\bm{\Gamma}_{\text{noise}}}^2,
\end{equation*}
where $\bm{\Gamma}_{\text{noise}}$ denotes the covariance matrix of the noise. Our objective is to infer the possible values of $u$ from the observation $\bm{d}$ within the Bayesian framework. To this end, we consider additive Gaussian noise $\bm{\epsilon}\sim\mathcal{N}(0,\bm{\Gamma}_{\text{noise}})$, where $\bm{\Gamma}_{\text{noise}}=\tau^{-1}\textbf{I}$ with $\tau^{-1}=(0.1\|\mathcal{S}\mathcal{G}u\|_{\infty})^2$. The reference Gaussian measure is $\mu_0=\mathcal{N}(0,\mathcal{C}_0)$, with covariance operator $\mathcal{C}_0=(\mathrm{I}-\alpha\Delta)^{-2}$ and $\alpha=1$, where $\Delta$ denotes the Laplace operator on $\Omega=[0,2\pi]$ with periodic boundary conditions.

We next specify the priors used in Bayes' formula. In this numerical example, we consider two distinct priors: the reference prior $\mu_0$ before the f-CNF transformation and the f-CNF prior $\mu_{f_{\theta}}$ obtained through direct data-driven training.

To train the f-CNF network, we construct a direct dataset $\mathbb{D}_{\mathbf{train}}^1$, similar to that in \cite{li2020fourier}, with the detailed definition given below:

\[
u(\theta) = \ln\left(\frac{\sin(3\theta) + 5}{3}\right) + \frac{\xi_0}{2\sqrt{2\pi}} + \sum_{k=1}^{10} \frac{\xi_k}{2(1+k^2)} \big(\sin(k\theta) + \cos(k\theta)\big),
\]
where $\xi_k$ ($k=0,1,\dots,10$) are independent standard normal random variables, i.e., $\xi_i \sim \mathcal{N}(0,1)$. We generate 2000 samples by sampling $\{\xi_k\}_{k=0}^{10}$ to form the training set $\mathbb{D}_{\mathrm{train}}$. The ground truth $u^{\dagger}$ is generated using the same procedure.

\begin{figure}[ht!]
	\centering
	\hspace{-0.03\textwidth} 
	\subfloat[Samples from $\mu_0$]{
		\includegraphics[ keepaspectratio=true, width=0.32\textwidth, clip=true]{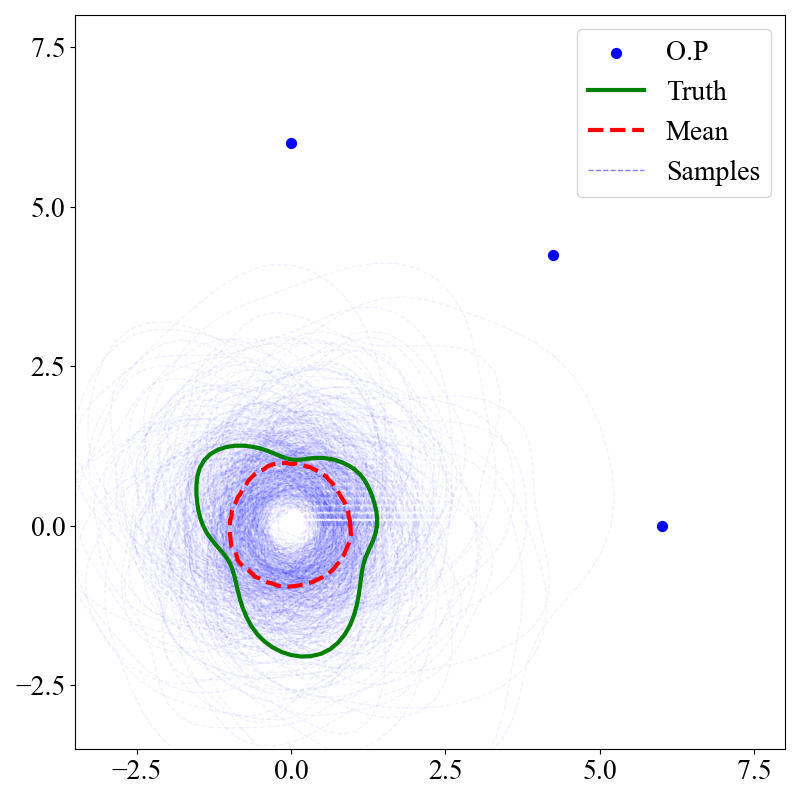}} 
  \hspace{0.08\textwidth} 
	\subfloat[Samples from $\mu_{f_{\theta}}$]{
		\includegraphics[ keepaspectratio=true, width=0.32\textwidth, clip=true]{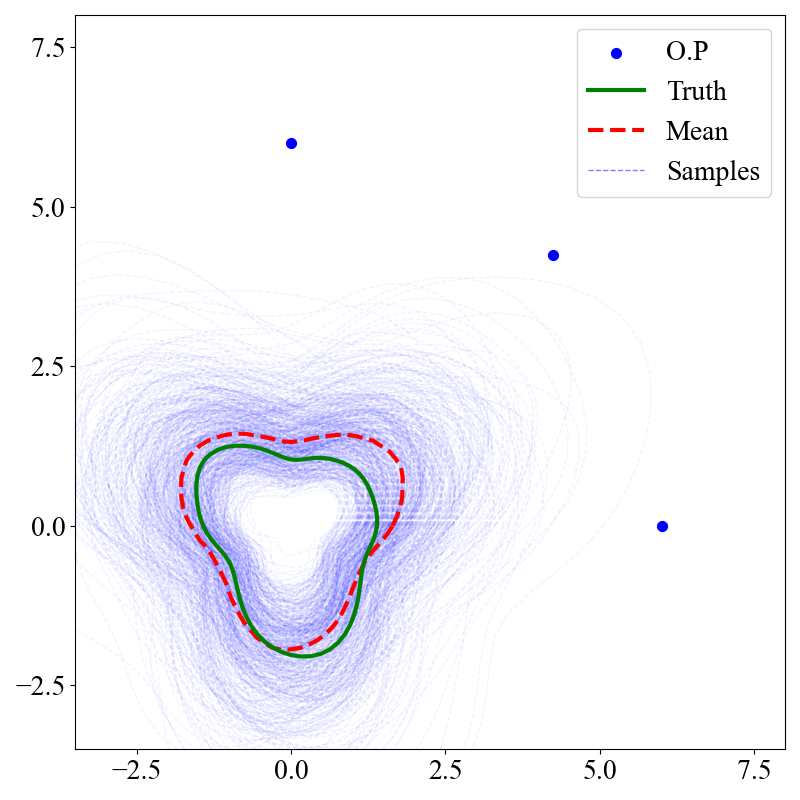}} 
	\caption{\emph{\small Samples from $\mu_0$ and trained f-CNF prior $\mu_{f_{\theta}}$. The green solid line represent the background truth $u^{\dagger}$, and the red dashed line represent the mean of the priors. The blue discrete points represent the position of the observation points (O.P) of our inverse problem. (a): Samples from $\mu_0$. (b): Samples from f-CNF prior $\mu_{f_{\theta}}$.}}
	
       \label{pic:hemholze prior}
\end{figure}

We employ Algorithm \ref{alg:generate model} to train the f-CNF prior $\mu_{f_{\theta}}$ on the training dataset $\mathbb{D}_{\textbf{train}}$ in infinite-dimensional space. Figure \ref{pic:hemholze prior} presents the ground truth $u^{\dagger}$, together with samples drawn from the initial prior $\mu_0$ and the learned prior $\mu_{f_{\theta}}$. It can be observed that, compared with those from $\mu_0$, the samples from $\mu_{f_{\theta}}$ exhibit greater similarity to the background truth $u^{\dagger}$. This indicates that $\mu_{f_{\theta}}$ provides more informative prior information for the posterior distribution.

\begin{figure}[ht!]
	\centering
	\hspace{-0.03\textwidth} 
	\subfloat[Posterior ($\mu_0$)]{
		\includegraphics[ keepaspectratio=true, width=0.32\textwidth, clip=true]{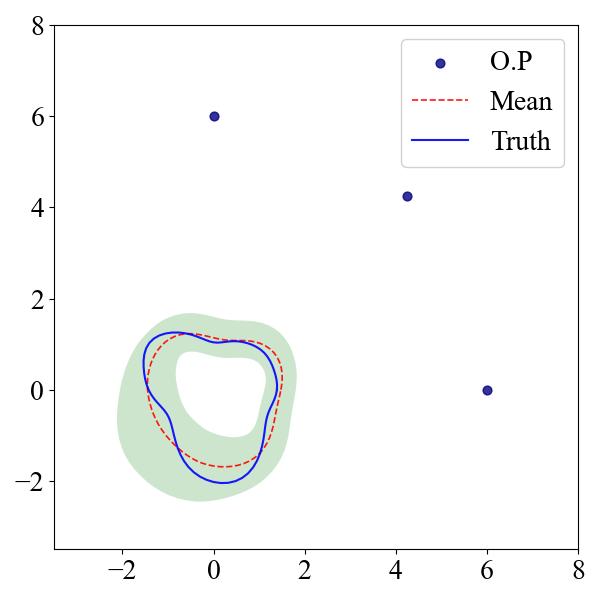}} 
  \hspace{0.08\textwidth} 
	\subfloat[Posterior ($\mu_{f_{\theta}}$)]{
		\includegraphics[ keepaspectratio=true, width=0.32\textwidth, clip=true]{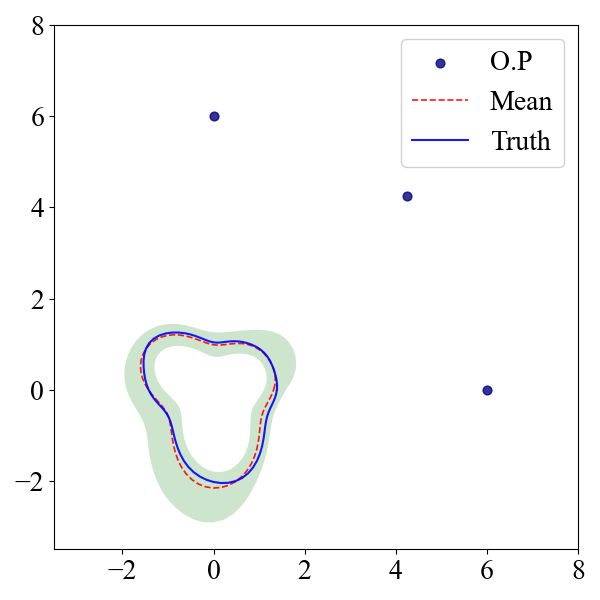}} 

  	\subfloat[Covariance ($\mu_0$)]{
		\includegraphics[ keepaspectratio=true, width=0.42\textwidth, clip=true]{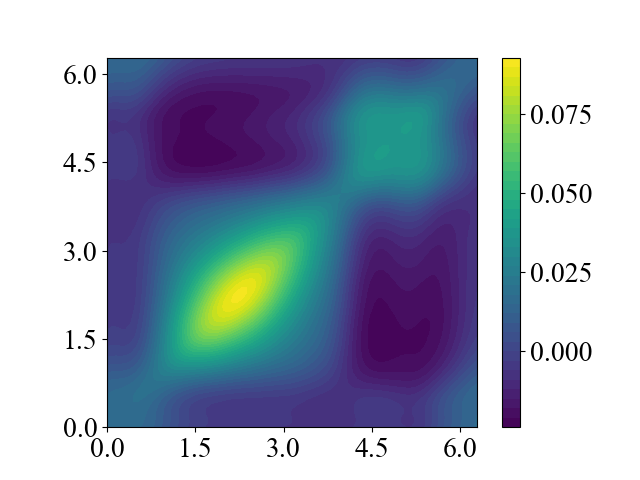}} 
	\subfloat[Covariance ($\mu_{f_{\theta}}$)]{
		\includegraphics[ keepaspectratio=true, width=0.42\textwidth, clip=true]{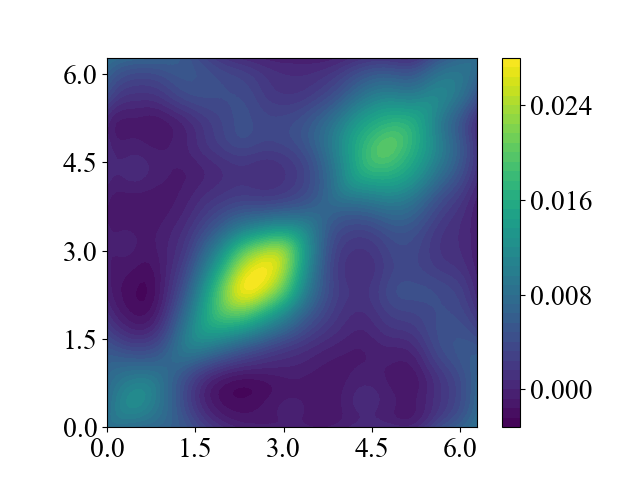}} 
	
	\caption{\emph{\small The comparison of the Bayesian posteriors with $\mu_0$, $\mu_{f_{\theta}}$ priors. (a)(b): The corresponding Bayesian posteriors with $\mu_0$ and $\mu_{f_{\theta}}$ priors. The blue solid line denotes the background truth $u^{\dagger}$, the red dashed line represents the mean of the posteriors. The green shaded region represents the $95\%$ credibility region of estimated posteriors, and the blue discrete points represent the position of the observation points (O.P). (c)(d): The covariance function of posteriors with $\mu_0$ and $\mu_{f_{\theta}}$ priors.}}
	
       \label{pic:hemholze posterior}
\end{figure}

We then draw posterior samples via the SMC-GM method (introduced in Subsection \ref{subsec:f-CNF prior}) under both the Gaussian prior $\mu_0$ and the trained prior $\mu_{f_\theta}$.
 For both methods, we generate 20000 samples. Numerical results are shown in Figure \ref{pic:hemholze posterior}. We also compute the $L^2$ error between the posterior mean and ground truth $u^\dagger$ for quantitative evaluation, as summarized in Table \ref{table:scar relative}.

\begin{table}[ht!]
	\renewcommand{\arraystretch}{1.5}
	\centering
		\caption{\emph{\small The relative errors between the background truth and the mean of posteriors.}} 
	\begin{tabular}{c|cc}
		\hline $\text{Relative Error}$ & posterior ($\mu_0$ prior)   & posterior ($\mu_{f_{\theta}}$ prior)  \\
		\hline $\text{Background Truth}$ $u^{\dagger}$& $0.17312$& $0.07151$\\
		\hline
	\end{tabular}
 \label{table:scar relative}
\end{table}
 
From Table \ref{table:scar relative} and Figure \ref{pic:hemholze posterior}, the posterior associated with $\mu_0$ yields accurate reconstructions and tight credible region near the three upper-right observation points, where the likelihood is informative. By contrast, in the unobserved lower-left region with limited likelihood information, reconstructions degrade and posterior uncertainty rises. The posterior derived from the learned f-CNF prior $\mu_{f_\theta}$ alleviates this limitation to some extent, yielding higher predictive accuracy and reduced uncertainty compared with the posterior of $\mu_0$.

\subsection{Inverse Heat Conduction Problem}\label{subsec:inverse_heat}

Inverse heat conduction problems have been studied for decades in scientific and engineering applications \cite{beck1985inverse}. In such problems, the temperature or heat flux on an inaccessible boundary is inferred from temperature measurements inside a solid body \cite{feng2018adaptive}. As our final numerical example, we consider the following heat equation:
$$
\frac{\partial w}{\partial t} = \frac{\partial}{\partial x} \left( h(x) \frac{\partial w}{\partial x} \right)
$$
with initial condition $w(x,0) = 0$. Here, $x$ and $t$ denote the spatial and temporal variables, respectively, $w(x, t)$ is the temperature, $h(x)$ is the thermal conductivity at each location, and $L$ is the length of the medium. All variables are given in dimensionless units. In this numerical experiment, the thermal conductivity $h(x)$ is defined as $h(x) = 0.5 + 0.15 \exp\left( -\left( \frac{x - 0.5}{0.2} \right)^2 \right).$

We assume that a heat flux is imposed on the left boundary $(x = 0)$ for $t \in [0,T_1]$, leading to the Neumann boundary condition:
$
\frac{\partial w}{\partial x}(0, t) = u_{\textbf{all}}(t),
$
where
$$
u_{\textbf{all}}(t)
=
\begin{cases}
u(t), & t \in [0,T_1], \\
0, & t > T_1.
\end{cases}
$$
Furthermore, we assume that the right end of the rod is insulated, which gives the boundary condition:
$
\frac{\partial w}{\partial x}(L, t) = 0.
$

Suppose that a temperature sensor is placed at the right end of the medium, i.e., $x=L$. The goal is to infer the heat flux $u(t)$ for $t \in [0,T_1]$ from the temperature measured by the sensor over the time interval $t\in [0,T_2]$. The corresponding PDE governing this process is given by
$$
\begin{cases}
\dfrac{\partial w}{\partial t}
=\dfrac{\partial}{\partial x} \left( h(x) \dfrac{\partial w}{\partial x} \right), \quad &\text{for } x \in [0, L], \ t \in [0, T_2], \\[7pt]
\dfrac{\partial w}{\partial x}(0, t) = u_{\textbf{all}}(t), \qquad &\text{for } \ t \in [0, T_2], \\[7pt]
\dfrac{\partial w}{\partial x}(L, t) = 0, \qquad &\text{for } \ t \in [0, T_2].
\end{cases}
$$

We now summarize the formulation of the above inverse heat conduction problem. Let $u(t)$ denote the model parameter, and let $\mathcal{G}: u(\cdot) \mapsto w(L,\cdot)$ denote the solution operator of the heat equation. The measurement operator $\mathcal{S}$ maps the solution $w(L,t)$ to its temporal observations at discrete instants $t\in[0,T_2]$, namely
$$\mathcal{S}: w(L,t) \mapsto \big(w(L,t_1),\,w(L,t_2),\,\dots,\,w(L,t_P)\big),$$
where $t_1,t_2,\dots,t_P \in [0,T_2]$ denote the measurement time instants. The noisy data vector $\bm{d}$ follows the forward mapping
\begin{align}\label{equ:heat inverse}
\bm{d} = \mathcal{S}\mathcal{G}(u) + \bm{\epsilon},
\end{align}
where $\bm{\epsilon}$ corresponds to additive Gaussian measurement noise.

For the inverse problem \eqref{equ:heat inverse}, the heat flux $u(t)$ takes time to propagate along the rod from left to right. Consequently, to accurately recover $u(t)$ over the entire interval $[0,T_1]$, it is typically necessary to take $T_2>T_1$ and observe $w(t)$ over the full interval $[0,T_2]$. Here, we choose $T_1 = 1$ and $T_2 = 2$. Consistent with Subsections \ref{subsec:simplesmooth} and \ref{subsec:inverse_scar}, we solve this problem within the Bayesian framework on infinite-dimensional spaces, using $\mu_0 = \mathcal{N}(0, \mathcal{C}_0)$ as the reference prior, where $\mathcal{C}_0 = (\text{I} - \alpha \Delta)^{-2}$ and $\alpha = 0.05$ is a prescribed constant. The Laplace operator is defined on $\Omega$ with homogeneous Neumann boundary conditions.

In this numerical experiment, we assume access to a large indirect dataset. Our objective is to train an f-CNF prior $\mu_{f_{\theta}}$ suitable for the inverse problem using Algorithm \ref{alg:generate model indirect}, and to demonstrate the effectiveness of the proposed method through comparison with the baseline prior $\mu_0$. We construct the indirect training dataset $\mathbb{D}_{\textbf{train}}^2$ similar to \cite{li2020fourier, akyildiz2025efficient}, with their detailed definitions given below.
\begin{itemize}
    \item Sample the parameter functions according to 
   \begin{equation*}
   u(t) = \cos(3\pi t) + 0.1\sum_{k=1}^{20} \frac{\xi_k \cos(k\pi t)}{1+4k^2},\quad \xi_k \sim \mathcal{N}(0,1).
    \end{equation*} The set $\{u_i\}_{i=1}^{10000}$ constitutes the parameter training dataset $\mathbb{D}^1$ generated by this sampling strategy.
    \item Define the indirect training dataset as $\mathbb{D}_{\textbf{train}}^2=\{\bm{d}_i \mid \bm{d}_i=\mathcal{S}_{\textbf{data}}\mathcal{G}(u_i)+\bm{\epsilon}_{\textbf{data}},\, u_i \in \mathbb{D}^1\}$. The measurement points of the operator $\mathcal{S}_{\textbf{data}}$ are taken as $\{ w(t_i) \mid i=1,2,\dots,50 \}$ with $t_i = \frac{2}{49}(i-1)$.
    Here, $\bm{\epsilon}_{\textbf{data}}$ denotes Gaussian white noise with fixed variance $0.001$.
\end{itemize}
The background truth $u^{\dagger}(t)$ is generated using the same sampling procedure as that used to construct $\mathbb{D}^1$.

\begin{figure}[ht!]
	\centering
	\hspace{-0.03\textwidth} 
	\subfloat[Samples from $\mu_0$]{
		\includegraphics[ keepaspectratio=true, width=0.40\textwidth, clip=true]{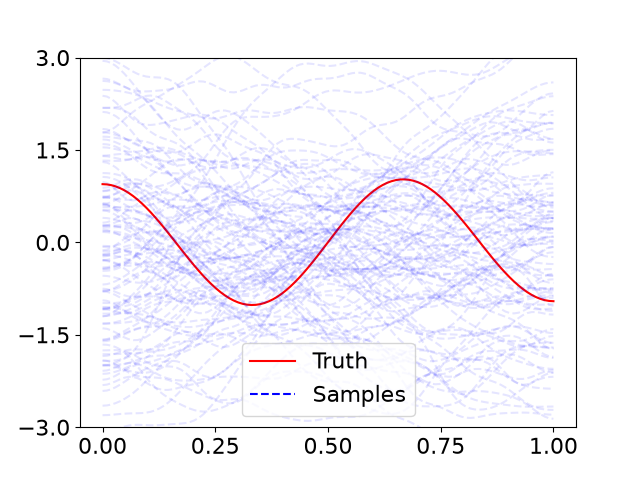}} 
  \hspace{0.08\textwidth} 
	\subfloat[Samples from $\mu_{f_{\theta}}$]{
		\includegraphics[ keepaspectratio=true, width=0.40\textwidth, clip=true]{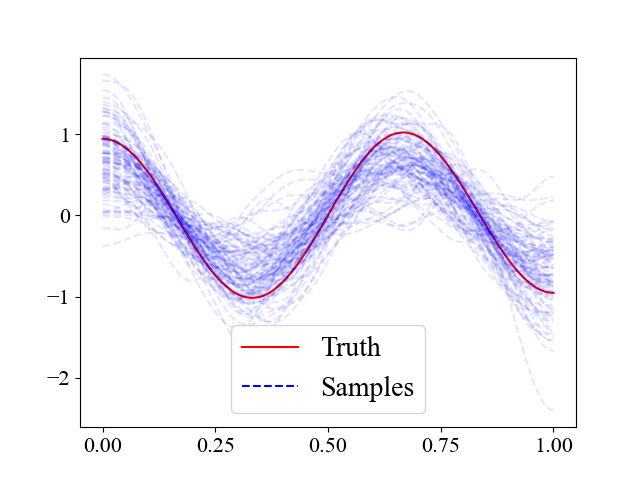}} 
	\caption{\emph{\small Samples from $\mu_0$ and trained f-CNF prior $\mu_{f_{\theta}}$. The red solid line represent the background truth of our inverse problem. (a): Samples from $\mu_0$. (b): Samples from f-CNF prior $\mu_{f_{\theta}}$.}}
	
       \label{pic:heat prior}
\end{figure}
Using $\mathbb{D}_{\textbf{train}}^2$, we train an f-CNF prior $\mu_{f_{\theta}}$ via Algorithm \ref{alg:generate model indirect}. To visualize the learned prior, Figure \ref{pic:heat prior} presents samples drawn from $\mu_0$ and $\mu_{f_{\theta}}$, together with their comparison to the background truth $u^{\dagger}$. The samples from $\mu_{f_{\theta}}$ exhibit greater similarity to $u^{\dagger}$ than those from $\mu_0$.

\begin{figure}[ht!]
	\centering
	
	\subfloat[Posterior (full)]{
		\includegraphics[ keepaspectratio=true, width=0.31\textwidth, clip=true]{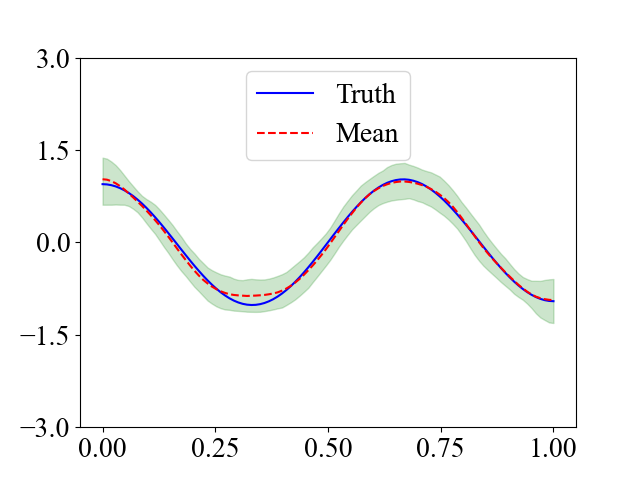}} 
	\subfloat[Posterior ($\mu_0$)]{
		\includegraphics[ keepaspectratio=true, width=0.31\textwidth, clip=true]{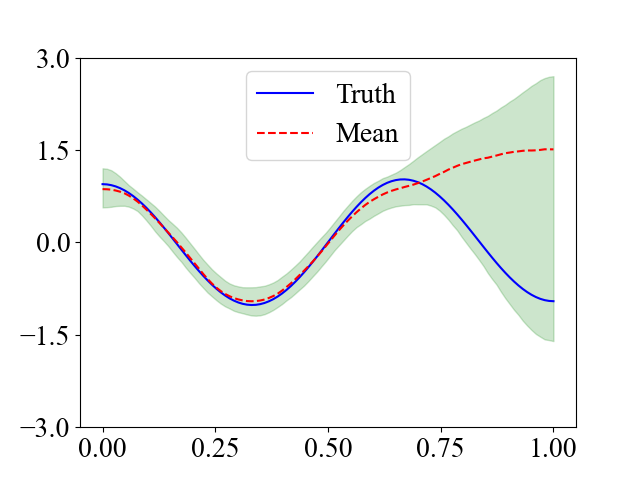}} 
        \subfloat[Posterior ($\mu_{f_{\theta}}$)]{
		\includegraphics[ keepaspectratio=true, width=0.31\textwidth, clip=true]{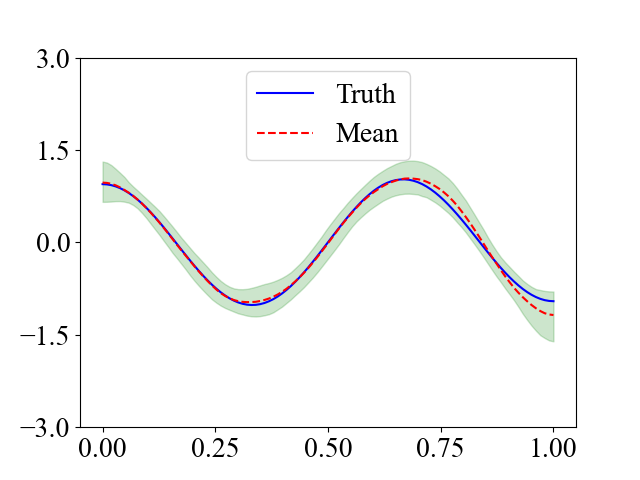}} 
	\\
	\subfloat[Covariance (full)]{
		\includegraphics[ keepaspectratio=true, width=0.31\textwidth, clip=true]{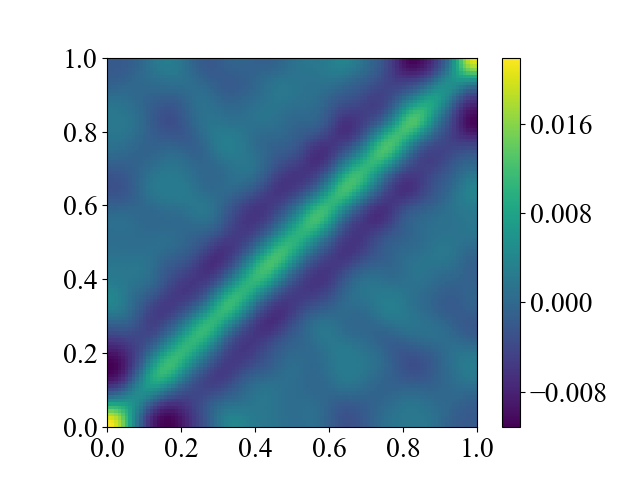}} 
	\subfloat[Covariance ($\mu_0$)]{
		\includegraphics[ keepaspectratio=true, width=0.31\textwidth, clip=true]{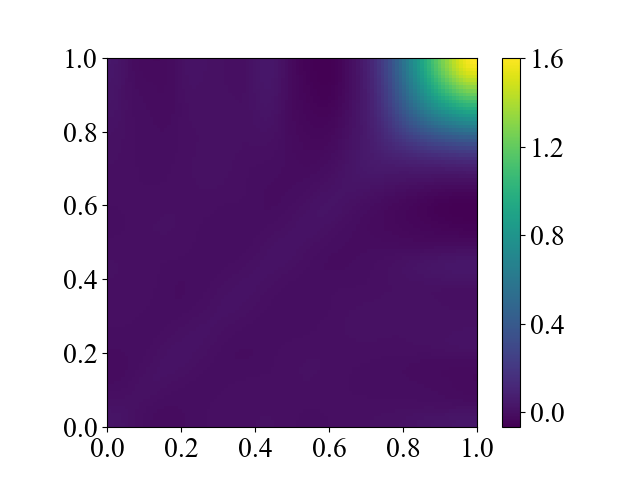}} 
        \subfloat[Covariance ($\mu_{f_{\theta}}$)]{
		\includegraphics[ keepaspectratio=true, width=0.31\textwidth, clip=true]{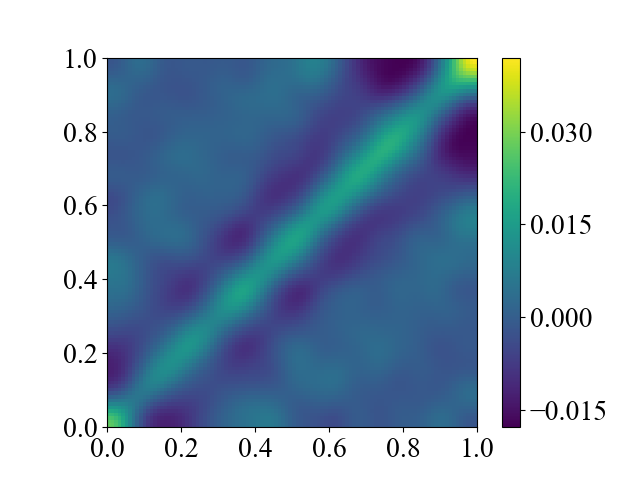}} 
	\caption{\emph{\small  Comparison of Bayesian posteriors: full measurements with $\mu_0$ prior, sparse measurements with $\mu_0$, $\mu_{f_{\theta}}$ priors (a): The Bayesian posterior with full measurement and $\mu_0$ prior. The blue solid line denotes the background truth $u^{\dagger}$, the red dashed line represents the mean of the posteriors. The green shaded region represents the $95\%$ credibility region of estimated posteriors. (b)(c): The corresponding Bayesian posteriors of sparse measurements with $\mu_0$, $\mu_{f_{\theta}}$ priors.
 (d)(e)(f): Posterior covariance functions under full measurements with $\mu_0$ prior and sparse measurements with $\mu_0$ and $\mu_{f_{\theta}}$ priors, respectively. }}
       \label{pic:heat posterior}
\end{figure}
We next perform Bayesian inversion using both $\mu_0$ and the learned f-CNF prior $\mu_{f_{\theta}}$. To further demonstrate the effectiveness of the trained prior $\mu_{f_{\theta}}$, we conduct three numerical experiments.

In the first numerical experiment, the measurement operator fully observes the solution $w(L,t)$ for $t\in[0,T_2]$, denoted by $\mathcal{S}_1: w(L,t) \mapsto (w(L,t^1), \dots, w(L,t^{50}))$, where $t^i =\frac{2}{49}(i-1)$ for $i=1,2,\dots,50$. Meanwhile, we assume that the measurements are contaminated by 20\% Gaussian random noise $\bm{\epsilon}_1 \sim \mathcal{N} (0,\bm{\Gamma}_{\text{noise}}^1)$, where $\bm{\Gamma}_{\text{noise}}^1 = \tau^{-1}\mathbf{I}$ and $\tau^{-1} = (0.20 \lVert \mathcal{S}_1\mathcal{G}u\rVert_{\infty})^2$. The measurement process is given by
$$
\bm{d}=\mathcal{S}_1\mathcal{G}(u)+\bm{\epsilon}_1.
$$
We choose $\mu_0$ as the prior and solve the inverse problem using SMC-GM. In the SMC-GM algorithm, 20,000 samples are generated, and the results are shown in Figure \ref{pic:heat posterior}. It can be observed that, when sufficient observational information is available, the prior $\mu_0$ yields accurate inversion results.

In the second and third numerical experiments, we shorten the measurement time interval. Specifically, we assume that the measurement operator observes the solution $w(L,t)$ only for $t\in[0,T_3]$, where $T_3=0.9$. This implies that no observational information is available after $t=0.9$, and thus the values of $u(t)$ for $t\in[0.9,1]$ cannot be accurately estimated from the measurements alone. We denote the corresponding measurement operator by $\mathcal{S}_2: w(L,t) \mapsto \big(w(L,t^1), \dots, w(L,t^{20})\big)$ where $t^i = 0.1 + \frac{8}{190}(i-1)$. Meanwhile, we assume that the measurements are contaminated by 20\% Gaussian random noise $\bm{\epsilon}_2 \sim \mathcal{N}(0, \bm{\Gamma}_{\text{noise}}^2)$, where $\bm{\Gamma}_{\text{noise}}^2 = \tau^{-1}\mathbf{I}$ and $\tau^{-1} = (0.20 \lVert \mathcal{S}_2\mathcal{G}u\rVert_{\infty})^2$.

In the second and third numerical experiments, we use $\mu_0$ and $\mu_{f_{\theta}}$ as the priors, respectively, and solve the inverse problem using SMC-GM algorithm. For both priors, 20,000 samples are generated, and the results are presented in Figure \ref{pic:heat posterior}.

In addition to the graphical results, we also compute the $L^2$ errors between the posterior means obtained in the three experiments and the background truth $u^{\dagger}$. The results are reported in Table \ref{table:heat relative}.
\begin{table}[ht!]
	\renewcommand{\arraystretch}{1.5}
	\centering
		\caption{\emph{\small The relative errors between the background truth and the mean of posteriors.}} 
	\begin{tabular}{c|ccc}
		\hline $\text{Relative Error}$ & posterior (full)   & posterior ($\mu_0$) & posterior ($\mu_{f_{\theta}}$) \\
		\hline $\text{Background Truth}$ $u^{\dagger}$& $0.07799$& $1.06424$& $0.10167$ \\
		\hline
	\end{tabular}
 \label{table:heat relative}
\end{table}
Combining the results in Table \ref{table:heat relative} and Figure \ref{pic:heat posterior}, we observe that, in the time intervals not covered by the measurements, the posterior corresponding to the prior $\mu_0$ exhibits extremely high uncertainty, and its posterior mean also deviates significantly from the background truth $u^{\dagger}$. In contrast, the learned prior $\mu_{f_{\theta}}$ effectively compensates for this deficiency, thereby demonstrating the effectiveness of the proposed prior learning method.

\section{Conclusion}\label{sec:Conclusion}

This work develops an f-CNF framework for constructing data-driven priors in infinite-dimensional Bayesian inverse problems. In particular, we establish the equivalence of the f-CNF prior to a Gaussian reference measure and derive an explicit expression for the RN derivative, providing a rigorous foundation for prior learning and posterior inference. Building on these results, we develop a measure-theoretic training framework for directly observed functional data and adopt a Sinkhorn-distance-based strategy for learning priors from indirect observations. For Bayesian inference with the resulting f-CNF priors, we establish the well-posedness of the posterior and adapt pCN and SMC-GM for posterior sampling. Numerical experiments on three representative PDE inverse problems demonstrate the effectiveness and applicability of the proposed framework in settings involving either direct or indirect observations.

However, the computational cost of ODE-based training remains a practical limitation. Future research will focus on developing more expressive and computationally efficient architectures, extending the framework to more challenging inverse problems, and investigating measure equivalence and characterizations of RN derivatives for broader classes of infinite-dimensional generative models, including diffusion-based approaches \cite{diffusionf2}.

\appendix

\section{Measure Equivalence and
Radon--Nikodym Derivatives}\label{appendix:Measure Equivalence and
Radon--Nikodym Derivatives}

\subsection{Fredholm and Regularized Determinants}\label{subappendix:Fredholm and Regularized Determinants}

Let $\mathcal{H}$ be a separable Hilbert space and $\mathcal{K}$ a compact operator on $\mathcal{H}$. We recall two determinant definitions.

\textbf{Definition 1 (Fredholm–Carleman determinant)}
When $\mathcal{K}$ is a trace-class operator, $\det_{1}(I+\mathcal{K})$ is well-defined as the absolutely convergent infinite product over all eigenvalues $\{\lambda_n\}$ of $\mathcal{K}$ (counted with algebraic multiplicity):
\begin{align*}
    \det\nolimits_{1}(I+\mathcal{K})=\prod_{n=1}^{\infty}\big(1+\lambda_n\big).
\end{align*}

\textbf{Definition 2 (Regularized Fredholm–Carleman determinant)}
Let $\mathcal{K}$ be a Hilbert–Schmidt operator on the Hilbert space $\mathcal{H}$. Denote $\lambda_n$ as the eigenvalues of $\mathcal{K}$. The standard infinite product for $\det\nolimits_1(I+\mathcal{K})$ may fail to converge in this setting, so the regularized determinant is defined as:
\[
\det\nolimits_{2}(I+\mathcal{K})=\prod_{n=1}^{\infty}(1+\lambda_n)e^{-\lambda_n}.
\]

For comprehensive theoretical discussions and proofs, see Section 6.4 in \cite{gaussianmeasure}.

\subsection{Proof of Lemma~\ref{thm:discrete}}\label{subappendix:proof of fnf}

\begin{proof}
Since the parameter $\theta$ does not play an explicit role in the discussion, we omit explicit reference to it. Accordingly, we adopt the following abbreviations:
\begin{itemize}
\item The operator $f_{\theta_n}^{(n)}$ is abbreviated as $f^{(n)}$.
\item The operator $\mathcal{F}_{\theta_n}^{(n)}$ is abbreviated as $\mathcal{F}^{(n)}$.
\item The operators $f_{\theta}$ and $\mathcal{T}_{\theta}$ are abbreviated as $f$ and $\mathcal{T}$, respectively.
\end{itemize}

The proof proceeds in three steps. Step 1: Express the composite map as $I+\mathcal{T}$ and show that $\operatorname{Im}(\mathcal{T})\subset\mathcal{H}$. Step 2: Establish $\mu\sim\mu_0$ by verifying the conditions of the Gaussian measure transformation theorem. Step 3: Derive the explicit RN density.

\textbf{Step 1.}
Since $f^{(n)}=I+\mathcal{F}^{(n)}$, successive expansion gives
$
f^{(N)}\circ\cdots\circ f^{(1)}
=I+\mathcal{T},
$
where
$
\mathcal{T}
=\mathcal{F}^{(1)}
+\sum_{n=2}^{N}\mathcal{F}^{(n)}
\circ f^{(n-1)}\circ\cdots\circ f^{(1)}.
$
Since $\operatorname{Im}(\mathcal{F}^{(n)})
\subset\mathcal{H}$ for all $n$, it follows that
$\operatorname{Im}(\mathcal{T})\subset\mathcal{H}$.

\textbf{Step 2.}
By the chain rule, $I+\mathcal{T}$ is differentiable along $\mathcal H$, with
\[
D_{\mathcal H}(I+\mathcal{T})(u)
=\Gamma^{(N)}\circ\cdots\circ\Gamma^{(1)},
\]
where
$
\Gamma^{(n)}=I_{\mathcal H}+D_{\mathcal H}\mathcal{F}^{(n)}(u_{n-1}),
\quad
u_{n-1}=f^{(n-1)}\circ\cdots\circ f^{(1)}(u),
$
for $n=2,\ldots,N$, with $u_0=u$ and
$\Gamma^{(1)}=I_{\mathcal H}+D_{\mathcal H}\mathcal{F}^{(1)}(u)$.

For every $u\in\mathcal{H}_u$, the point spectrum of $D_{\mathcal H}\mathcal{F}^{(n)}(u)$ is disjoint from $(-\infty,-1]$, which implies that each $\Gamma^{(n)}$ is injective for all $n=1,\dots,N$. Consequently, $D_{\mathcal H}(I+\mathcal{T})(u)$ is invertible for all $u\in\mathcal{H}_u$.

Moreover, $f=I+\mathcal{T}$ is a homeomorphism, and $\mathcal{T}\in\mathcal{HC}^1(\mu_0,\mathcal H)$ with $D_{\mathcal H}\mathcal{T}(u)$ trace class. Applying Theorem 6.6.7 and Corollary 6.6.8 from \cite{gaussianmeasure}, we establish the mutual absolute continuity $\mu\sim\mu_0$, where $\mu=\bigl(f^{(N)}\circ f^{(N-1)}\circ\cdots\circ f^{(1)}\bigr)_{\#}\mu_0$.

\textbf{Step 3.} We divide the proof into two cases.

\noindent\textit{Case $N=1$.} Let $\mathcal F^{(1)}=\mathcal F$. By Theorem 6.6.7 and Corollary 6.6.8 in \cite{gaussianmeasure}, with the Gaussian divergence understood in the localized sense described after Lemma 5.8.10 of \cite{gaussianmeasure}, the Radon--Nikodym derivative is given by
\begin{align*}
\frac{d\bigl((I+\mathcal{F})_{\#} \mu_0 \bigr)}{d\mu_0}\bigl((I+\mathcal{F})(u)\bigr)
=\frac1{\Lambda_{\mathcal{F}}\left(u\right)},
\end{align*}
where
\begin{align*}
\Lambda_{\mathcal{F}}(u)
&:= \bigl|\det\nolimits_2\bigl(I_{\mathcal H}+D_{\mathcal H}\mathcal{F}(u)\bigr)\bigr|
\exp\left[\delta_{\mu_0}\mathcal{F}(u)-\frac12\lVert\mathcal{F}(u)\rVert_{\mathcal{H}}^2\right] \\
&= \bigl|\det\nolimits_1\bigl(I_{\mathcal H}+D_{\mathcal H}\mathcal{F}(u)\bigr)\bigr|
\exp\left[-\langle u,\mathcal{F}(u) \rangle_{\mathcal{H}}^{\mathrm{ext}}-\frac12\lVert\mathcal{F}(u)\rVert_{\mathcal{H}}^2\right].
\end{align*}

\noindent\textit{Case $N>1$.} Let $f = I+\mathcal{T}$ as defined in Subsection \ref{subsec:RN derivative}. The chain rule implies
\begin{align*}
D_{\mathcal H}f(u)
&= D_{\mathcal H}(f^{(N)} \circ \cdots \circ f^{(1)})(u) \\
&= D_{\mathcal H}f^{(N)}\bigl(f^{(N-1)}\circ \cdots \circ f^{(1)}(u)\bigr)
\circ D_{\mathcal H}f^{(N-1)}\bigl(f^{(N-2)}\circ\\
&\cdots \circ f^{(1)}(u)\bigr)
\cdots \circ D_{\mathcal H}f^{(1)}(u).
\end{align*}
Since the determinant of a composition factorizes, we have
\[
\det\nolimits_{1}\bigl(D_{\mathcal H}f(u)\bigr)
= \det\nolimits_{1}\bigl(D_{\mathcal H}f^{(N)}(u_{N-1})\bigr)\det\nolimits_{1}\bigl(D_{\mathcal H}f^{(N-1)}(u_{N-2})\bigr)
\cdots \det\nolimits_{1}\bigl(D_{\mathcal H}f^{(1)}(u)\bigr),
\]
where $u_k = f^{(k)}\circ\cdots \circ f^{(1)}(u)$ for $k=1,\dots,N-1$.
Recall that $\mathcal{T}(u)=f(u)-u$, applying Corollary 6.6.8 from \cite{gaussianmeasure} yields
\begin{align*}
\frac{d\bigl(f_{\#} \mu_0 \bigr)}{d\mu_0}\bigl(f(u)\bigr)
&=\frac{1}{\bigl|\det_{1}\bigl(I_{\mathcal H}+D_{\mathcal H}\mathcal{T}(u)\bigr)\bigr|}
\exp\left[\frac12\lVert\mathcal{T}(u)\rVert_{\mathcal{H}}^2 + \langle u,\mathcal{T}(u) \rangle_{\mathcal{H}}^{\mathrm{ext}} \right] \\
&=\prod_{k=1}^{N}\bigl|\det\nolimits_{1}(D_{\mathcal H}f^{(k)}(u_{k-1}))\bigr|^{-1}
\exp\left[\frac12\lVert\mathcal{T}(u)\rVert_{\mathcal{H}}^2 + \langle u,\mathcal{T}(u) \rangle_{\mathcal{H}}^{\mathrm{ext}} \right].
\end{align*}

This completes the proof of the lemma.

\end{proof}


\subsection{Proof of Theorem~\ref{thm:limits}}\label{subappendix:proof of limit}
\begin{proof}
Let $\nu$ be the signed measure defined by
\[
\nu(A) = \int_{A} \mathcal{R}(u) \, \mu_0(du).
\]
By the definition of the Radon–Nikodym derivative, we have
$
\frac{d\nu}{d\mu_0}(u) = \mathcal{R}(u).
$
Since $\mathcal{R}(u) > 0$ pointwise, the measure $\nu$ is equivalent to $\mu_0$.

Let $\phi \in C_b(\mathcal{H}_u)$. We compute
\begin{align*}
\lim_{N\to\infty}\int_{\mathcal{H}_u}\phi(u)\,d{f_N}_{\#}\mu_0
= \lim_{N\to\infty}\int_{\mathcal{H}_u}\phi(f_N(u))\,d\mu_0 
= \int_{\mathcal{H}_u}\phi(f(u))\,d\mu_0 
= \int_{\mathcal{H}_u}\phi(u)\,df_{\#}\mu_0.
\end{align*}
The interchange of the limit and the integral is justified by the dominated convergence theorem. This implies that ${f_N}_{\#} \mu_0$ converges weakly to $f_{\#} \mu_0$.

Next, we rewrite the limit using the density:
\begin{align*}
\lim_{N\to\infty}\int_{\mathcal{H}_u}\phi(u)\,d{f_N}_{\#}\mu_0
= \lim_{N\to\infty}\int_{\mathcal{H}_u}\phi(u)\mathcal{R}_N(u)\,d\mu_0 
= \int_{\mathcal{H}_u}\phi(u)\mathcal{R}(u)\,d\mu_0 
= \int_{\mathcal{H}_u}\phi(u)\,d\nu.
\end{align*}
Again, the dominated convergence theorem allows us to exchange the limit and the integral. Hence, ${f_N}_{\#}\mu_0$ converges weakly to $\nu$.

By the uniqueness of weak limits, we conclude that $f_{\#}\mu_0 = \nu$.
\end{proof}

\subsection{Proof of Theorems \ref{thm:RN derivative}}\label{subappendix:proof of RN}

\begin{proof}
For simplicity, we omit the explicit dependence on $\theta$. Let $\Delta t=T/N$ and $t_n=n\Delta t$. For any $u\in\mathcal H_u$, define the Euler approximation by
\[
u_N^0(u)=u,\qquad
u_N^{n+1}(u)=u_N^n(u)+\Delta t\,h(u_N^n(u),t_n),
\quad n=0,\ldots,N-1.
\]
Set $f^N(u)=u_N^N(u)$ and $\mathcal T^N(u)=f^N(u)-u$. By the continuity and global Lipschitz assumptions on $h$, the Euler approximation converges to the exact flow $f$, namely $\lim\limits_{N \rightarrow \infty} f^N(u) = f(u)$
 for every $u\in\mathcal H_u$.

Each Euler step is given by $f^{(n)}=I+\Delta t\,h(\cdot,t_n)$. For sufficiently large $N$, we have $\Delta t L<1$ and $\Delta t M_2<1$. The former ensures that $f^{(n)}$ is a homeomorphism by the contraction mapping theorem, while the latter guarantees the invertibility of $I_{\mathcal H}+\Delta t\,D_{\mathcal H}h(u,t_n)$. Moreover, the assumptions of $h$ imply that $\Delta t\,h(\cdot,t_n)$ is continuous in $\mathcal H$, and its  derivative is a trace class operator and trace-norm continuous. Hence Lemma \ref{thm:discrete} applies to $f^N=f^{(N-1)}\circ\cdots\circ f^{(0)}$, and $(f^N)_{\#}\mu_0\sim\mu_0$ for all sufficiently large $N$.

Since $\mathcal T^N(u)=\Delta t\sum_{n=0}^{N-1}h(u_N^n(u),t_n)$, the assumptions imply
$\|\mathcal C_0^{-1}\mathcal T^N(u)\|_{\mathcal H_u}\le TM_1$ and
$\|D_{\mathcal H}\mathcal T^N(u)\|_1\le(1+\Delta t M_2)^N$.
Thus, the inner product yields
\[
\langle u,\mathcal T^N(u)\rangle_{\mathcal H}^{\mathrm{ext}}
=\langle u,\mathcal C_0^{-1}\mathcal T^N(u)\rangle_{\mathcal H_u}.
\]

For any $y\in\mathcal H_u$, let
$u_N^*=(f^N)^{-1}(y)$ and $u^*=f^{-1}(y)$. We introduce
\[
A_n^N=D_{\mathcal H}h(u_N^n(u_N^*),t_n),
\qquad z_N=\mathcal T^N(u_N^*).
\]
Applying Lemma \ref{thm:discrete}, the RN density of the discrete pushforward is
\begin{align*}
\mathcal R_N(y)
:=\frac{d(f^N_{\#}\mu_0)}{d\mu_0}(y)=&\exp\Biggl(
-\sum_{n=0}^{N-1}\log|\det\nolimits_1(I_{\mathcal H}+\Delta t A_n^N)|
+\frac12\|z_N\|_{\mathcal H}^2\\
&\qquad+\langle u_N^*,\mathcal C_0^{-1}z_N\rangle_{\mathcal H_u}
\Biggr).
\end{align*}

We next establish the pointwise convergence of $\mathcal R_N$. Denote by $v(t;u^*)$ the exact solution with initial condition $v(0;u^*)=u^*$. Since each Euler step has an inverse Lipschitz constant bounded by $(1-\Delta t L)^{-1}$, we have
\[
\lim\limits_{N \rightarrow \infty}\|u_N^*-u^*\|_{\mathcal H_u}
\le\lim\limits_{N \rightarrow \infty}(1-\Delta t L)^{-N}
\|f^N(u^*)-f(u^*)\|_{\mathcal H_u}
=0.
\]
Consequently, the stability and convergence of the Euler scheme imply
\[
\lim\limits_{N \rightarrow \infty}\max_{0\le n\le N}
\|u_N^n(u_N^*)-v(t_n;u^*)\|_{\mathcal H_u}
=0.
\]

For sufficiently small trace-class $A$, the Fredholm determinant satisfies
$\log|\det_1(I_{\mathcal H}+A)|=\operatorname{Tr}_{\mathcal H}A+O(\|A\|_1^2)$.
Using the trace-norm continuity of $D_{\mathcal H}h$ and its uniform bound $M_2$, we obtain
\[
\lim\limits_{N \rightarrow \infty}\sum_{n=0}^{N-1}
\log|\det\nolimits_{1}(I_{\mathcal H}+\Delta t A_n^N)|
=
\int_0^T\operatorname{Tr}_{\mathcal H}
D_{\mathcal H}h(v(t;u^*),t)\,dt.
\]

Since $\lim\limits_{N \rightarrow \infty}z_N=y-u^*=\mathcal T(u^*)$ in $\mathcal H_u$ and
$\|\mathcal C_0^{-1}z_N\|_{\mathcal H_u},
\|\mathcal C_0^{-1}\mathcal T(u^*)\|_{\mathcal H_u}\le TM_1$, we have
$
\lim\limits_{N \rightarrow \infty}\|z_N-\mathcal T(u^*)\|_{\mathcal H}^2
\le\lim\limits_{N \rightarrow \infty}2TM_1\|u_N^*-u^*\|_{\mathcal H_u}
=0.
$
Moreover, the boundedness of $\mathcal C_0^{-1}z_N$, together with
$z_N\to\mathcal T(u^*)$ in $\mathcal H_u$ and the injectivity of
$\mathcal C_0$, yields
$\mathcal C_0^{-1}z_N\rightharpoonup
\mathcal C_0^{-1}\mathcal T(u^*)$ weakly in $\mathcal H_u$.
Combining this with the strong convergence $u_N^*\to u^*$ in
$\mathcal H_u$, we obtain
\[
\lim\limits_{N \rightarrow \infty}
\langle u_N^*,\mathcal C_0^{-1}z_N\rangle_{\mathcal H_u}
=
\langle u^*,\mathcal C_0^{-1}\mathcal T(u^*)\rangle_{\mathcal H_u}.
\]
Combining these limits, we obtain
$\lim\limits_{N\to\infty}\mathcal R_N(y)=\mathcal R(y)$, where
\begin{align*}
\mathcal R(y)
&=\exp\Biggl(
-\int_0^T\operatorname{Tr}_{\mathcal H}
D_{\mathcal H}h(v(t;u^*),t)\,dt+\frac12\|\mathcal T(u^*)\|_{\mathcal H}^2
+\langle u^*,\mathcal C_0^{-1}\mathcal T(u^*)\rangle_{\mathcal H_u}
\Biggr).
\end{align*}

It remains to construct an integrable dominating function. Let
$M_3=\|\mathcal C_0^{1/2}\|M_1$ and
$M_4=\|\mathcal C_0\|M_1$.
The assumptions imply
\[
\|z_N\|_{\mathcal H}\le TM_3,\qquad
\|z_N\|_{\mathcal H_u}\le TM_4,\qquad
\|u_N^*\|_{\mathcal H_u}\le\|y\|_{\mathcal H_u}+TM_4.
\]
Furthermore, the logarithmic determinant expansion and the uniform trace-norm bound yield, for all sufficiently large $N$,
\[
\left|\sum_{n=0}^{N-1}
\log|\det\nolimits_{1}(I_{\mathcal H}+\Delta t A_n^N)|\right|
\le TM_2+C,
\]
where $C$ is independent of $N$ and $y$. Consequently,
\begin{align*}
0<\mathcal R_N(y)
&\le\exp\Bigl(
TM_2+C+\tfrac12T^2M_3^2+TM_1(\|y\|_{\mathcal H_u}+TM_4)
\Bigr)\\
&\le C_1\exp(C_2\|y\|_{\mathcal H_u}).
\end{align*}
By Fernique's theorem, $g(y)=C_1\exp(C_2\|y\|_{\mathcal H_u})$ belongs to $L^1(\mu_0)$. Therefore, applying Theorem \ref{thm:limits} to the Euler approximations yields
$\frac{d(f_{\#}\mu_0)}{d\mu_0}(y)=\mathcal R(y)$. Since $\mathcal R(y)>0$ for $\mu_0$-almost every $y$, we conclude that $f_{\#}\mu_0\sim\mu_0$. Finally, setting $y=f(u)$ gives $u^*=u$, and using
$
\langle u,\mathcal T(u)\rangle_{\mathcal H}^{\mathrm{ext}}
=\langle u,\mathcal C_0^{-1}\mathcal T(u)\rangle_{\mathcal H_u}
$
yields the desired RN derivative.
\end{proof}

\section{Posterior Sampling and Computational Details}\label{app:sampling}\label{appendix:Posterior Sampling and Computational Details}
\subsection{pCN Sampling under f-CNF Priors}\label{subappendix:pCN Sampling under f-CNF Priors}
Here, we briefly introduce the pCN algorithm for Bayesian posterior inference with the f-CNF prior. For a detailed presentation of the algorithm, we refer readers to \cite{cotter2013,dashti2013bayesian}.

Recall the original Bayesian formula associated with the f-CNF prior:
\begin{align}\label{equ:bayes_f_CNF}
\frac{d\mu}{d\mu_{f_{\theta}}}(u) = \frac{1}{Z_{\mu}^{\theta}} \exp(-\Phi(\bm{d}|u)).
\end{align}
Here, $Z_{\mu}^{\theta}$ denotes a finite positive constant defined as
$
Z_{\mu}^{\theta} = \int_{\mathcal{H}_u} \exp (-\Phi(\bm{d}|u))\mu_{f_{\theta}}(du).
$
Additionally, the Radon–Nikodym derivative of the measure $\mu_{f_{\theta}}$ with respect to the Gaussian measure $\mu_0$ is given by
$
\frac{d\mu_{f_{\theta}}}{d\mu_0}(u) = \exp(\mathcal{W}_{\theta}(u)).
$

Since the pCN algorithm requires the prior measure to be Gaussian, we reformulate the original Bayesian formula as
\begin{align}\label{equ:bayesnewprior}
\frac{d\mu}{d\mu_{0}}(u) = \frac{1}{Z_{\mu}^{\theta}} \exp\left(-\Phi(\bm{d}|u)+\mathcal{W}_{\theta}(u)\right).
\end{align}
The application of the pCN algorithm to this reformulated posterior measure is summarized in Algorithm \ref{alg:CNF-pCN}.

\begin{algorithm}[htbp]
  \caption{f-CNF pCN}
  \label{alg:CNF-pCN}
  \begin{algorithmic}[1]
    \REQUIRE 
    f-CNF prior $\mu_{f_{\theta}}$ (transformed from Gaussian $\mu_0$),
    RN potential $\mathcal{W}(u)$ with $\frac{d\mu_{f_{\theta}}}{d\mu_0}(u) = \exp\left(\mathcal{W}(u)\right)$,
    likelihood potential $\Phi(\bm{d}|u)$,
    initial state $u_0$,
    total samples $K$,
    step size $\beta$
    \ENSURE 
    Posterior samples $\{u_k\}_{k=1}^K$.

    \STATE \textbf{Initialization}: Set $k = 0$, start from initial state $u_0$
    \REPEAT
        
        \STATE \textbf{Sample base noise}: $v_0 \sim \mu_0$ (Gaussian base measure)
        \STATE \textbf{Proposal generation}: $v = \sqrt{1 - \beta^2} \, u_k + \beta \, v_0$
        \STATE \textbf{Acceptance probability}:
        \[
        a(u_k, v) = \min\left\{
        \exp\big(
        \Phi(\bm{d}|u_k) - \Phi(\bm{d}|v)
        + \mathcal{W}(v) - \mathcal{W}(u_k)
        \big),\; 1
        \right\}
        \]
        \STATE \textbf{Sample update}:
        $u_{k+1} =
        \begin{cases}
        v & \text{with probability } a(u_k, v), \\
        u_k & \text{with probability } 1-a(u_k, v)
        \end{cases}$
        \STATE $k = k + 1$
    \UNTIL{$k = K$}

    \RETURN $\{u_k\}_{k=1}^K$
  \end{algorithmic}
\end{algorithm}

\subsection{SMC-GM algorithm}\label{subappendix:SMC-GM algorithm}
The pCN algorithm efficiently samples from posterior distributions. However, when the observational data are highly informative and the posterior differs significantly from the prior, a small step size is required to maintain the algorithm's acceptance rate. This leads to strong correlations between successive samples, making it challenging to obtain independent samples. To resolve this issue, the SMC-GM algorithm is proposed in \cite{lu2025sequential}. The fundamental algorithmic principles are established in \cite{lu2025sequential}, and the complete implementation procedure for Bayesian posterior inference under the f-CNF prior is summarized in Algorithm \ref{alg:smc_gm}. Here, SMC-GM provides an approximate posterior sampler, with its approximation properties supported by the results in \cite{lu2025sequential}. The same work also describes a Metropolis--Hastings correction that uses the fitted Gaussian mixture as a proposal and preserves the tempered posterior at each mutation step \cite{lu2025sequential}. In this paper, our numerical study focuses on the effect of the learned f-CNF prior on Bayesian inversion. For computational efficiency, we use the SMC-GM approximation, which draws directly from the fitted mixture without a Metropolis--Hastings acceptance step.

\begin{algorithm}[htbp]
  \caption{f-CNF SMC-GM}
  \label{alg:smc_gm}
  \begin{algorithmic}[1]
    \REQUIRE 
    f-CNF prior $\mu_{f_{\theta}}$, 
    number of particles $N$, 
    step sizes $\{h_j\}_{j=1}^J$ (sum to 1), 
    potential function $\Phi(\bm{d}|u)$
    \ENSURE 
    Particles $\{u_J^{(i)}\}_{i=1}^N$ approximating posterior $\mu$ of our Bayes formula.

    \STATE \textbf{Initialization}: 
    Sample $u_0^{(i)} \sim \mu_{f_{\theta}}$, $i=1,\dots,N$

    \FOR{$j = 0$ \TO $J-1$}
        \STATE Set weights $W_j^{(i)} = 1/N$, $i=1,\dots,N$

        \STATE \textbf{Reweighting}:
        $
        \widetilde{W}_{j+1}^{(i)} = W_j^{(i)} \exp\left( -h_{j+1} \,\Phi\left(u_j^{(i)}\right) \right)
        $
        \\
        Normalize: $W_{j+1}^{(i)} = \widetilde{W}_{j+1}^{(i)} \big/ \sum_{l=1}^N \widetilde{W}_{j+1}^{(l)}$

        \STATE \textbf{Resampling}:
        Resample $\widehat{u}^{(i)}$ from $\{u_j^{(l)}\}$ with weights $\{W_{j+1}^{(l)}\}$

        \STATE \textbf{Gaussian Mixture Fitting}:
        Fit $\widehat{\mu}_j = \sum_{k=1}^M \pi_k \mathcal{N}(m_k, \mathcal{C}_k)$ to $\{\widehat{u}^{(i)}\}$

        \STATE \textbf{Mutation (GM Sampling)}:
        Sample $u_{j+1}^{(i)} \sim \widehat{\mu}_j$
    \ENDFOR

    \RETURN $\{u_J^{(i)}\}_{i=1}^N$
  \end{algorithmic}
\end{algorithm}

\subsection{Accelerate}\label{subappendix:Accelerate}
When evaluating the RN derivative for a sample $u$, the required quantity in our algorithm is $\exp\left(\mathcal{W}_{\theta}(u)\right) = \frac{d\mu_{f_{\theta}}}{d\mu_0}(u)$. However, the direct formula from Theorem \ref{thm:RN derivative} evaluates this derivative at the forward-transformed state, yielding $\exp(\mathcal{W}_{\theta}(f_{\theta}(u))) = \frac{d\mu_{f_{\theta}}}{d\mu_0}(f_{\theta}(u))$ rather than at $u$.

A conceptually straightforward but computationally inefficient approach to obtain $\mathcal{W}_{\theta}(u)$ using Theorem \ref{thm:RN derivative} is to first determine the pre-image of $u$ under the flow map $f_{\theta}$. Let $v = f_{\theta}^{-1}(u)$ such that $f_{\theta}(v) = u$. Applying Theorem \ref{thm:RN derivative} to this pre-image $v$ then gives $\mathcal{W}_{\theta}(f_{\theta}(v)) = \mathcal{W}_{\theta}(u)$, which is the desired quantity.

Although mathematically valid, this two-step procedure is computationally prohibitive. Specifically, it requires solving the inverse ODE to obtain $v = f_{\theta}^{-1}(u)$, followed by solving the forward ODE to evaluate the trajectory needed for Theorem \ref{thm:RN derivative}. This necessitates two complete ODE solutions per sample $u$. In the context of Markov chain Monte Carlo (MCMC) sampling, where such evaluations must be performed iteratively for numerous samples, this overhead results in significant computational inefficiency. To address this bottleneck, we present Theorem \ref{thm:generate model}. This theorem provides an alternative formulation that evaluates
$\mathcal{W}_\theta(u)$ in a single backward ODE solve,
avoiding the separate inverse then forward procedure.

\begin{theorem}\label{thm:generate model}
Consider the forward mapping $f_{\theta}$ as defined in Theorem \ref{thm:RN derivative}. For any $u \in \mathcal{H}_u$, let $v(t)$ denote the solution of the following ODE:
\begin{align}\label{equ:thm21}
\left\{
\begin{array}{l}
\displaystyle \frac{dv(t)}{dt} = h(v(t),t; \theta), \\
v(T) = u.
\end{array}
\right.
\end{align}
and define $\mathcal{T}_{\theta}(u) = f_{\theta}(u)-u$. Then, the RN derivative of the pushforward measure $\mu_{f_{\theta}}$ with respect to $\mu_0$ is given by:
$$\frac{d\mu_{f_{\theta}}}{d\mu_0}(u)=\exp\left(\frac{1}{2}\lVert\mathcal{T}_{\theta}(v(0))\rVert_{\mathcal{H}}^2 + \langle v(0),\mathcal{T}_{\theta}(v(0)) \rangle_{\mathcal{H}}^{\mathrm{ext}} - \int_0^{T} \mathrm{Tr}_{\mathcal{H}}(D_{\mathcal{H}}h(v(t),t; \theta)) dt\right).$$
\end{theorem}

\noindent\textbf{Proof of Theorem \ref{thm:generate model}}

\begin{proof}
Since the solution $v(t)$ of the underlying ODE exists and is unique, it follows that
$
f_{\theta}^{-1}(u) = v(0).
$
Substituting this into the RN derivative yields
\begin{align*}
\frac{d\mu_{f_{\theta}}}{d\mu_0}(u) &= \frac{d\mu_{f_{\theta}}}{d\mu_0}(f_{\theta}(f_{\theta}^{-1}(u))) \\
&= \exp\left(\frac{1}{2}\lVert\mathcal{T}_{\theta}(v(0))\rVert_{\mathcal{H}}^2 + \langle v(0), \mathcal{T}_{\theta}(v(0)) \rangle_{\mathcal{H}}^{\mathrm{ext}} - \int_0^{T} \text{Tr}_{\mathcal{H}} \, D_{\mathcal{H}}h(v(t), t; \theta) \, dt\right),
\end{align*}
where $v(t)$ denotes the solution of the corresponding ODE \eqref{equ:thm21}.
\end{proof}

It can be observed that computing the RN derivative according to the original Theorem \ref{thm:RN derivative} necessitates solving the following ODE system:
\begin{align*}
    \left\{
    \begin{array}{l}
    \displaystyle \frac{dv(t)}{dt} = h(v(t),t; \theta), \\
    v(0) = f_{\theta}^{-1}(u).
    \end{array}
    \right.
\end{align*}
This approach requires solving an additional inverse ODE to obtain $f_{\theta}^{-1}(u)$. In contrast, the proposed Theorem \ref{thm:generate model} only requires solving
\begin{align}\label{equ:ODE2}
\left\{
\begin{array}{l}
\displaystyle \frac{dv(t)}{dt} = h(v(t),t; \theta), \\
v(T) = u.
\end{array}
\right.
\end{align}
This bypasses the need to compute $f_\theta^{-1}(u)$
in a separate preliminary step. Since
$\mathcal T_\theta(v(0))=u-v(0)$, the RN derivative
can be evaluated using the same backward trajectory,
without an additional forward ODE solve.
Consequently, only a single ODE solve is required,
thereby reducing the overall computational overhead.

\subsection{Splitting acceleration strategy}\label{subappendix:Splitting acceleration strategy}

The proposed f-CNF prior shares a similar measure-theoretic structure with the TV-Gaussian prior \cite{tvgaussian}, both admitting a density with respect to a reference Gaussian measure $\mu_0$:
\[
\frac{d\mu_{\mathrm{prior}}}{d\mu_0}(u)
=\exp(\mathcal{W}_{\theta}(u)),
\]
where $\mathcal{W}_{\theta}$ denotes the corresponding log-density. For the TV-Gaussian prior, it is determined by the total variation of $u$, whereas for the f-CNF prior, it is characterized by the learned flow.

This structure affects the efficiency of pCN sampling. With the reference Gaussian prior, the acceptance probability depends only on the negative log-likelihood $\Phi(u)$. In contrast, for the f-CNF prior, it also depends on $\mathcal{W}_{\theta}(u)$, which may vary significantly between proposed states. Consequently, large proposal steps can lead to high rejection rates and slow mixing. Motivated by this similarity to the TV-Gaussian prior, the splitting acceleration strategy (S-pCN) \cite{tvgaussian} can be adapted to improve sampling efficiency under the f-CNF prior.

\noindent\textbf{Detailed procedure of the splitting acceleration algorithm}

Let $u_{\text{current}}$ denote the current state. The algorithm proceeds in two consecutive stages: ``local updates'' and a ``global acceptance test''.

\textbf{Stage 1: Multiple local updates based on $\mathcal{W}_{\theta}(u)$}

Initialize the intermediate state $v_0 = u_{\text{current}}$, and perform $k$ successive pCN updates as follows.
\begin{enumerate}
    \item Generate a proposal using the standard pCN proposal rule:
    \[
    v_{\text{prop}} = \sqrt{1-\beta^2}\,v_{i-1} + \beta w,
    \]
    where $w\sim\mu_0$ and $\beta$ denotes the step-size parameter.
    \item Compute the acceptance probability based solely on the energy functional:
    \[
    \mathrm{acc}_{\mathcal{W}_{\theta}}\big(v_{\text{prop}},v_{i-1}\big) = \min\Big\{1,\exp\big(\mathcal{W}_{\theta}(v_{\text{prop}})-\mathcal{W}_{\theta}(v_{i-1})\big)\Big\}.
    \]
    \item Update the intermediate state according to the computed acceptance probability to obtain the $i$-th iterate $v_i$.
\end{enumerate}
After $k$ iterations, we obtain the intermediate state $v_k$. This stage requires only evaluations of the low-cost energy functional $\mathcal{W}_{\theta}(u)$, which permits a large number of local iterations $k$ to be performed efficiently.

\textbf{Stage 2: Final global acceptance test}

Treat $v_k$ as the proposed state and compute the final acceptance probability using the data-fidelity term $\Phi(u)$:
\[
\mathrm{acc}_\Phi\big(v_k,u_{\text{current}}\big) = \min\Big\{1,\exp\big(-\Phi(v_k)+\Phi(u_{\text{current}})\big)\Big\}.
\]
Perform the final accept–reject step using this probability to produce the updated sample $u_{\text{next}}$.

This splitting strategy decomposes each sampling iteration into the two stages described above. The first stage performs numerous low-cost local updates and accept–reject decisions relying entirely on the energy term $\mathcal{W}_{\theta}(u)$. The second stage then performs a single global acceptance test that incorporates the data-fidelity term $\Phi(u)$. The method mitigates the slow mixing induced by the small step-size constraint in standard pCN algorithms and improves overall sampling efficiency, all without introducing any approximation error.

\section{Well-Posedness of the Bayesian Posterior}\label{appendix:Well-Posedness of the Bayesian Posterior}
\subsection{Proof of Theorem \ref{thm:wellposedness}}\label{subappendix:proof of wellposedness}

\begin{proof}
    Consider the corresponding Bayesian posterior distribution:
    \begin{align*}
        \frac{d\mu}{d\mu_{f_{\theta}}}(u)=\frac{1}{Z_{\theta}}\mathbb{L}(\bm{d}|u),
    \end{align*}
    where
    $
    \frac{d\mu_{f_{\theta}}}{d\mu_0}(u)=\exp(\mathcal{W}_{\theta}(u)).
    $
    By the chain rule, we have
    $$
    \frac{d\mu}{d\mu_0}(u)=\frac{1}{Z_{\theta}}\mathbb{L}(\bm{d}|u)\exp(\mathcal{W}_{\theta}(u)).
    $$
    Based on the proof of Theorem \ref{thm:RN derivative}, there exists a function $g(u)=C_1\exp(C_2||u||_{\mathcal{H}_u})$ such that
    $
    \exp(\mathcal{W}_{\theta}(u))\leq g(u).
    $
    Consequently, since the likelihood is bounded by a constant $M$, we have
    \begin{align*}
        \int_{\mathcal{H}_u}|\mathbb{L}(\bm{d}|u)\exp(\mathcal{W}_{\theta}(u))|\mu_0(du)\leq \int_{\mathcal{H}_u}M\times C_1\exp(C_2||u||_{\mathcal{H}_u})\mu_0(du)<\infty,
    \end{align*}
    which implies that $\mathbb{L}(\bm{d}|u)\exp(\mathcal{W}_{\theta}(u))\in L^1(\mathcal{H}_u,\mu_0)$, and that there exists a function $Mg(u)\in L^1(\mathcal{H}_u,\mu_0)$ such that, for any $\bm{d}$,
    $|\mathbb{L}(\bm{d}|u)\exp(\mathcal{W}_{\theta}(u))|\leq Mg(u).$
    Therefore, based on Theorem 3.6 of \cite{MR4624340}, we conclude that our Bayesian inverse problem is well-posed in the weak, Hellinger, and total variation metrics.

Next, observe that
$
||u||_{\mathcal{H}_u}^p \mathbb{L}(\bm{d}|u) \exp(\mathcal{W}_{\theta}(u))\leq ||u||_{\mathcal{H}_u}^p C_3 \exp(C_2 ||u||_{\mathcal{H}_u}).
$
For a sufficiently small $\epsilon > 0$, we have $\exp(\epsilon||u||^2)\in L^1(\mathcal{H}_u,\mu_0)$. Furthermore, for $||u||_{\mathcal{H}_u}$ sufficiently large, the quadratic growth dominates, yielding
$$
||u||_{\mathcal{H}_u}^p C_3 \exp(C_2 ||u||_{\mathcal{H}_u})<\exp(\epsilon ||u||_{\mathcal{H}_u}^2).
$$
Consequently, one can construct a function $h\in L^1(\mathcal{H}_u,\mu_0)$ such that
$$
||u||_{\mathcal{H}_u}^p \mathbb{L}(\bm{d}|u)\exp(\mathcal{W}_{\theta}(u))\leq h(u).
$$
Thus, based on Theorem 3.12 of \cite{MR4624340}, we conclude that our Bayesian inverse problem is well-posed with respect to the $p$-Wasserstein metric.
\end{proof}

\section{Functional Flow Architectures}\label{appendix:Functional Flow Architectures}

\subsection{Verification of the Continuous Planar Flow Conditions}\label{subappendix:Verification of the Continuous Planar Flow Conditions}

The functional continuous planar flow is defined by the following equation:
\begin{equation*}
\left\{
\begin{aligned}
\frac{dv(t)}{dt} &= \sum_{i=1}^{L}u_i(t)\sigma(\langle v(t),w_i(t)\rangle_{\mathcal H_u}+b_i(t)),\\
v(0)&=u_0.
\end{aligned}
\right.
\end{equation*}
Therefore, the vector field in Theorem \ref{thm:RN derivative} is given by
\[
h(v,t;\theta)=\sum_{i=1}^{L}u_i(t)\sigma(\langle v,w_i(t)\rangle_{\mathcal H_u}+b_i(t)),
\]
where $u_i(t)$, $w_i(t)$, and $b_i(t)$ are parameterized by neural networks, and $\theta$ collectively denotes the corresponding trainable network parameters. For brevity, denote
\[
\xi_i(v,t)=\langle v,w_i(t)\rangle_{\mathcal H_u}+b_i(t).
\]
Our primary objective is to prove that $h(\cdot,\cdot;\theta)\in\mathcal T(\mu_0,[0,T],\mathcal H)$. The verification consists of four parts.

\textbf{Part 1.} First, we verify that $h(\cdot,t;\theta)\in\mathcal{HC}^{1}(\mu_0,\mathcal H)$.
Recall that $u_i(t)$ and $w_i(t)$ are parameterized as
$
u_i(t)=\sum_{j=1}^{M}\alpha_i^{(j)}(t)\phi_j,
$
and
$
w_i(t)=\sum_{j=1}^{M}\beta_i^{(j)}(t)\phi_j,
$
where $\phi_j$ are eigenfunctions of $\mathcal C_0$ with positive eigenvalues. Consequently, $\phi_j\in\mathcal H$, and for any $t\in[0,T]$, $\theta\in\Theta$, and $v\in\mathcal H_u$, we have $h(v,t;\theta)\in\mathcal H$.

For any $g\in\mathcal H$, the Cameron--Martin derivative is
\[
D_{\mathcal H}h(v,t;\theta)g
=\sum_{i=1}^{L}\sigma'(\xi_i(v,t))\,u_i(t)
\langle g,w_i(t)\rangle_{\mathcal H_u}.
\]
Since
$
\langle g,w_i(t)\rangle_{\mathcal H_u}
=\langle g,\mathcal C_0w_i(t)\rangle_{\mathcal H},
$
we can write
\[
D_{\mathcal H}h(v,t;\theta)
=\sum_{i=1}^{L}\sigma'(\xi_i(v,t))\,
u_i(t)\otimes w_i(t),
\]
 Thus, $D_{\mathcal H}h(v,t;\theta)$ is a finite-rank Hilbert--Schmidt operator. Moreover, its continuity in $(v,t)$ with respect to the Hilbert--Schmidt norm follows from the continuity of $\sigma'$, $u_i(t)$, $w_i(t)$, and $b_i(t)$. Therefore, $h(\cdot,t;\theta)\in\mathcal{HC}^1(\mu_0,\mathcal H)$ for each $t\in[0,T]$. The same parameter continuity also ensures that $h$ is continuous on $\mathcal H_u\times[0,T]$.

\textbf{Part 2.} Second, we establish the boundedness condition. We have
\begin{align*}
\mathcal C_0^{-1}h(v,t;\theta)
&=\sum_{i=1}^{L}\mathcal C_0^{-1}u_i(t)\sigma(\xi_i(v,t))=\sum_{i=1}^{L}\sum_{j=1}^{M}
\alpha_i^{(j)}(t)\mathcal C_0^{-1}\phi_j\,\sigma(\xi_i(v,t)).
\end{align*}
Since $\alpha_i^{(j)}(t)$ is continuous on the compact interval $[0,T]$, the constants
$
C_1=\max_{i,j}\sup_{t\in[0,T]}|\alpha_i^{(j)}(t)|, C_2=\max_j\|\mathcal C_0^{-1}\phi_j\|_{\mathcal H_u}
$
are finite. Since $\sigma=\tanh$ satisfies $|\sigma|\le1$, we obtain
$
\|\mathcal C_0^{-1}h(v,t;\theta)\|_{\mathcal H_u}
\le LMC_1C_2=:M_1<\infty.
$

\textbf{Part 3.} We next bound the trace norm on the Cameron--Martin space $\mathcal H$. Writing
$
D_{\mathcal H}h(v,t;\theta)=\sum_{i=1}^{L}A_i$, where $A_i=\sigma'(\xi_i(v,t))\,
u_i(t)\otimes_{\mathcal H}\mathcal C_0w_i(t),
$
the rank-one property and $|\sigma'|\le1$ yield
$
\|A_i\|_1
=|\sigma'(\xi_i(v,t))|\,
\|u_i(t)\|_{\mathcal H}\,
\|\mathcal C_0^{1/2}w_i(t)\|_{\mathcal H_u}
\le M_i,
$
where $M_i<\infty$ exists by continuity on the compact interval $[0,T]$. Here we used $\|\mathcal C_0w_i(t)\|_{\mathcal H}=\|\mathcal C_0^{1/2}w_i(t)\|_{\mathcal H_u}$. Consequently,
$
\|D_{\mathcal H}h(v,t;\theta)\|_1
\le\sum_{i=1}^{L}\|A_i\|_1
\le\sum_{i=1}^{L}M_i=:M_2<\infty.
$
Furthermore, since each $A_i$ depends continuously on $(v,t)$ in the trace norm, so does the operator $D_{\mathcal H}h(v,t;\theta)$.

\textbf{Part 4.} Finally, we verify the global Lipschitz condition on $\mathcal H_u$. Since $|\sigma'|\le1$, for any $v_1,v_2\in\mathcal H_u$, we have
\begin{align*}
&\|h(v_1,t;\theta)-h(v_2,t;\theta)\|_{\mathcal H_u}\le\sum_{i=1}^{L}\|u_i(t)\|_{\mathcal H_u}
\bigl|\sigma(\xi_i(v_1,t))-\sigma(\xi_i(v_2,t))\bigr|\\
&\quad\le\left(\sum_{i=1}^{L}\|u_i(t)\|_{\mathcal H_u}\|w_i(t)\|_{\mathcal H_u}\right)
\|v_1-v_2\|_{\mathcal H_u}.
\end{align*}
Since $u_i(t)$ and $w_i(t)$ are continuous on $[0,T]$, the constant
\[
L_h=\sup_{t\in[0,T]}\sum_{i=1}^{L}
\|u_i(t)\|_{\mathcal H_u}\|w_i(t)\|_{\mathcal H_u}<\infty
\]
is finite. Hence
$
\|h(v_1,t;\theta)-h(v_2,t;\theta)\|_{\mathcal H_u}
\le L_h\|v_1-v_2\|_{\mathcal H_u}.
$

In summary, we have shown that $h(\cdot,\cdot;\theta)\in\mathcal T(\mu_0,[0,T],\mathcal H)$.

\subsection{Other flow model}\label{subappendix:Other flow model}

Although the numerical experiments employ the functional continuous planar flow, other flow models inspired by \cite{fnf} may also satisfy the conditions of Theorem \ref{thm:RN derivative}. Several such models are presented below, with detailed proofs omitted as they follow arguments analogous to those for the functional continuous planar flow.


\textbf{Functional continuous radial-basis flow:}
Let $\{\phi_j\}_{j=1}^{M}$ be the first $M$ eigenfunctions of
$\mathcal C_0$, chosen orthonormally in $\mathcal H_u$ and
corresponding to positive eigenvalues. Define the projection
and reconstruction operators by
\[
\mathcal P u
=\bigl(\langle u,\phi_1\rangle_{\mathcal H_u},
\ldots,\langle u,\phi_M\rangle_{\mathcal H_u}\bigr)^{\mathrm T},
\qquad
\mathcal Q z=\sum_{j=1}^{M}z_j\phi_j,
\]
for $u\in\mathcal H_u$ and $z\in\mathbb R^M$.
The functional continuous radial-basis flow is defined by
\[
\left\{
\begin{aligned}
\frac{dv(t)}{dt}
&=\mathcal Q\sum_{k=1}^{K}a_k(t)
\rho\!\left(\gamma_k(t)
\|\mathcal Pv(t)-c_k(t)\|_{\mathbb R^M}^{2}\right),\\
v(0)&=u_0,
\end{aligned}
\right.
\]
where $K$ is the number of radial basis functions,
$a_k(t),c_k(t)\in\mathbb R^M$ are coefficient vectors and
centers, and $\gamma_k(t)>0$ are scale parameters.
These parameters are continuous on $[0,T]$ and jointly
parameterized by a neural network with trainable parameters $\theta$.
The function $\rho(r)=\exp(-r)$.

\textbf{Functional continuous Sylvester flow:} 
Building upon the Sylvester flow in \cite{sylvesterflow} and the functional Sylvester flow from \cite{fnf}, we define the functional continuous Sylvester flow as follows:
\begin{align*}
\left\{
\begin{array}{l}
\displaystyle \frac{dv(t)}{dt} = \mathcal{A}(t)h(\mathcal{B}(t)v(t)+b(t)), \\
v(0) = u_0,
\end{array}
\right.
\end{align*}
where, for any fixed $t \in [0,T]$, $\mathcal{B}(t)$ is a bounded linear operator mapping $\mathcal{H}_u$ to $\mathbb{R}^M$, $\mathcal{A}(t)$ is a bounded linear operator mapping $\mathbb{R}^M$ to $\mathcal{H}_u$, $b(t) \in \mathbb{R}^M$, and $h(x)=\tanh(x)$. The time-dependent operators $\mathcal{A}(t)$, $\mathcal{B}(t)$, and bias $b(t)$ are parameterized by a neural network with trainable parameters $\theta$. We can parameterize
$\mathcal A(t)=\mathcal Q A_M(t)$ and
$\mathcal B(t)=B_M(t)\mathcal P$,
where $A_M(t),B_M(t)\in\mathbb R^{M\times M}$ and $b(t)$
depend continuously on $t\in[0,T]$.
Here, $\mathcal P$ and $\mathcal Q$ are the covariance-eigenbasis
projection and reconstruction operators.

\section*{Acknowledgments}
This research was partially supported by the National Natural Science Foundation of China (Grant Nos. 12322116 and 12271428) and by the Fundamental and Interdisciplinary Disciplines Breakthrough Plan of the Ministry of Education of China (Grant No. JYB2025XDXM101).

\bibliographystyle{amsplain}
\bibliography{references}

\end{document}